\RequirePackage{fix-cm}
\documentclass[10pt]{amsart}
\usepackage[utf8]{inputenc}
\usepackage[T1]{fontenc}
\usepackage{soul}

\usepackage{graphicx}
\usepackage{longtable}
\usepackage{float}
\usepackage{wrapfig}
\usepackage{textcomp}
\usepackage{marvosym}
\usepackage{latexsym}
\usepackage{amssymb}
\usepackage{xargs}
\usepackage{natbib}
\usepackage{algorithm2e}
\usepackage{xcolor}
\usepackage{subcaption}
\definecolor{HalfGray}{gray}{0.55}
\definecolor{OliveGreen}{rgb}{0,.35,0}
\definecolor{webbrown}{rgb}{.6,0,0}
\definecolor{BrightViolet}{rgb}{0.5,0.2,0.8}
\definecolor{Maroon}{cmyk}{0, 0.87, 0.68, 0.32}
\definecolor{RoyalBlue}{cmyk}{1, 0.50, 0, 0.25}
\definecolor{Black}{cmyk}{0, 0, 0, 0}

\definecolor{ccccccc}{RGB}{204,204,204}
\definecolor{c808080}{RGB}{128,128,128}
\definecolor{c999999}{RGB}{153,153,153}
\definecolor{ce6e6e6}{RGB}{230,230,230}

\usepackage[unicode]{hyperref}
\hypersetup{
colorlinks=true,
linktocpage=true,
pdfstartview=FitH,
breaklinks=true,
pdfpagemode=UseNone,
pageanchor=true,
pdfpagemode=UseOutlines,
plainpages=false,
bookmarksnumbered,
bookmarksopen=false,
bookmarksopenlevel=1,
hypertexnames=true,
pdfhighlight=/O,
urlcolor=Maroon, linkcolor=RoyalBlue, citecolor=OliveGreen, 
pdftitle={},
pdfauthor={},
pdfsubject={},
pdfkeywords={},
pdfcreator={pdfLaTeX},}
\usepackage{pgf}

\usepackage{pgfplots}
\usepackage{amsthm,amsmath,enumerate,bbm,tikz,tabularx,multicol,environ}
\usetikzlibrary{matrix,arrows,calc,patterns}
\usepgfplotslibrary{groupplots,dateplot}
\usetikzlibrary{patterns,shapes.arrows}
\pgfplotsset{compat=newest}

\newtheorem{theorem}{Theorem}[section]
\newtheorem{proposition}[theorem]{Proposition}
\newtheorem{corollary}[theorem]{Corollary}
\newtheorem{lemma}[theorem]{Lemma}
\theoremstyle{definition}
\newtheorem{definition}[theorem]{Definition}
\newcounter{assump}
\theoremstyle{remark}

\DeclareMathOperator*{\argmin}{arg\,min}

\author{Joon Kwon}
\date{\today}
\title{A regret minimization approach to fixed-point iterations}
\begin{document}

\begin{abstract}
  In Euclidean spaces, we propose a conversion scheme that turns regret minimizing
  algorithms into fixed point iterations, with convergence guarantees
  following from regret bounds. The connection builds upon the
  equivalence between fixed points of nonexpansive maps and zeros of
  cocoercive operators. The classical Krasnoselskii--Mann iterations
  are recovered by converting the Online Gradient Descent algorithm.
  We then convert algorithms from the AdaGrad family of regret
  minimizers, and obtain best-iterate residual guarantees with bounds
  that are adaptive to the best-in-hindsight scalar scaling (resp.\ diagonal/full-matrix preconditioning).


\end{abstract}
\maketitle

\section{Introduction}
\label{sec:introduction}

In optimization and numerical analysis, iterative methods play a
central role for solving complex and large-scale problems. Instead of
aiming for an exact solution, they output a sequence of points that
hopefully converges to a solution. They offer important computational advantages: cheap iteration cost, low
memory requirements, and great scalability. They have been successfully
used for solving linear systems~\citep{jacobi1845ueber,seidel1874uber}, ODEs and PDEs~\citep{runge1895numerische,kutta1901beitrag,richardson1910approximate,douglas1956numerical,young1954iterative},
optimization problems~\citep{cauchy1847methode,kuhn1951nonlinear}, and
in particular, have been an indispensable tool for the current AI
revolution by offering excellent scalability for neural network
optimization~\citep{robbins1951stochastic,mcmahan2010adaptive,duchi2011adaptive,tieleman2012lecture,kingma2015adam}.

\subsection{Fixed-point iterations}
\label{sec:fixed-point-iter}

Fixed points are a very general framework where for a given operator
$F$, a solution is defined as a point $x_*$ such that $F(x_*)=x_*$. In
this context, iterative methods are called \emph{fixed-point iterations}.
The most favorable case is when the operator
$F:\mathcal{X}\to \mathcal{X}$ is a \emph{contraction} in some complete metric space
$(\mathcal{X},\nu)$, meaning it is $L$-Lipschitz continuous for some coefficient
$0\leqslant L<1$. Then, the Banach--Picard theorem~\citep{banach1922operations}
guarantees that there exists a unique fixed point $x_*\in \mathcal{X}$ and that from
any initial point $x_1\in \mathcal{X}$, iterating operator $F$, setting
$x_{t+1}=F(x_t)$ for all $t\geqslant 1$, ensures that
$\nu(x_t,x_*)\leqslant L^{t-1} \nu(x_1,x_*)$ for all $t\geqslant 1$, which corresponds to
a geometric convergence to the fixed point. Applications include the proof of the Picard--Lindelöf
theorem~\citep{lindelof1894application}, computation of stationary
distributions of Markov chains~\citep{kemeny1960finite} and value
iteration in dynamic programming~\citep{bellman1957dynamic}.

An important relaxation is the case where the operator is only
$1$-Lipschitz continuous, which is also called \emph{nonexpansive}.
This situation is not as straightforward because a fixed point may not exist
(e.g.,\ for a translation), and even if a fixed point exists,
iterating the operator may not converge to the fixed point (e.g.,\ for
a rotation, if the initial point is not already the center of the
rotation). But remarkably, on a convex set in a Euclidean
space\footnote{Proposition~\ref{prop:km-intro} was first presented in
  the Hilbert setting by
   \citet{vaisman2005convergencia} and formally
  published in \citet{cominetti2014rate}. Prior to that, \citet{krasnosel1955two,schaefer1957uber,edelstein1966remark,ishikawa1976fixed,edelstein1978nonexpansive} established strong convergence
  in various kinds of normed spaces; weak convergence in uniformly
  convex Banach spaces
  with Fréchet differentiable norms was studied in
  \citet{reich1979weak}. In more general spaces,
having a vanishing fixed-point residual is 
a weaker property, called \emph{asymptotic regularity}, that has been
established under various assumptions by
\citet{goebel1983iteration,groetsch1972note}, and in general normed
spaces by \citet{baillon1996rate} and \citet{cominetti2014rate}.
For KM iterates of a nonexpansive self-map on a closed convex set with a
fixed point, asymptotic regularity implies weak convergence in a reflexive
Banach space satisfying Opial's property~\citep{opial1967weak}. Reflexivity
provides weak subsequential compactness of the bounded iterates.},
as soon as a fixed point exists, the
Krasnoselskii--Mann (KM)
iteration \citep{mann1953mean,krasnosel1955two}
\[ x_{t+1}=\frac{x_t+F(x_t)}{2},\quad t\geqslant 1, \]
which averages the previous point with its image by the operator,
guarantees convergence to a fixed point, as well as an upper bound
on the last-iterate \emph{fixed-point residual} that vanishes as $1/\sqrt{T}$,
where $T$ is the number of iterations.
\begin{proposition}[Krasnoselskii--Mann iterations \citep{vaisman2005convergencia,cominetti2014rate}]
\label{prop:km-intro}
Let $d\geqslant 1$ be an integer, $\mathcal{X}\subset \mathbb{R}^d$ a closed convex set, $F:\mathcal{X}\to \mathcal{X}$ a
$1$-Lipschitz continuous function for $\left\|\,\cdot\,\right\|_2$, $x_1\in \mathcal{X}$, and for $t\geqslant 1$,
\[ x_{t+1}=\frac{x_t+F(x_t)}{2}. \]
We assume that $F$ admits a fixed point $x_*\in \mathcal{X}$. Then, $(x_t)_{t\geqslant 1}$ converges to a fixed point of $F$ and for all $T\geqslant 1$,
\[ \left\| F(x_T)-x_T \right\|_2\leqslant \frac{2\left\| x_1-x_* \right\|_2}{\sqrt{T}}. \]
\end{proposition}
Note that it is not possible to deduce an upper bound on the distance of the
iterates to their limit from the above upper bound on the fixed-point residual
$\left\| F(x_T)-x_T \right\|_2$. The latter is a weaker measure, of how
far a point is from \emph{being} a fixed point.

As first noted by \citet{combettes2004solving}, many important
optimization algorithms are instances of KM iterations: Gradient
Descent~\citep{cauchy1847methode}, the Proximal Point
method~\citep{martinet1970breve}, the Forward-Backward
splitting~\citep{lions1979splitting}, the Douglas--Rachford
splitting~\citep{douglas1956numerical}, the Chambolle--Pock
algorithm~\citep{chambolle2011first}, etc., which can all be seen
through the lens of the very general \emph{monotone operator
  splitting} framework, which applies KM iterations to a wide range of
problems
\citep{bauschke2011convex,ryu2016primer,ryu2022large,combettes2024geometry}.

Remarkably, it has been recently proved that \emph{Halpern iterations}
with carefully chosen viscosity
parameters~\citep{lieder2021convergence} as well as KM iterations with
Nesterov's momentum~\citep{boct2023fast} guarantee vanishing
fixed-point residual at the faster $1/T$ rate for nonexpansive operators.

\subsection{The wish for \emph{adaptive} fixed-point iterations}
\label{sec:quest-adaptive-fixed}
Let $d\geqslant 1$ be an integer, $I$ the identity operator on
$\mathbb{R}^d$ and $F:\mathbb{R}^d\to \mathbb{R}^d$. For simplicity, we restrict this
discussion to operators defined on $\mathbb{R}^d$. An interesting remark is
that operators $F$ and $F_{L}:=I+\frac{1}{L}(F-I)$, for any $L\neq 0$, have the same fixed
points. We now restrict to $L>0$. Operator $F_{L}$ is then obtained from $F$ by a kind of
positive scaling that preserves the fixed points. In the case where
$F$ is not nonexpansive, there may exist a value $L>0$ such that $F_L$
is nonexpansive, and using the latter for performing KM iterations
would then be beneficial. Even in the case where $F$ \emph{is}
nonexpansive, performing KM iterations with $F_L$ instead of $F$ may
provide stronger guarantees and faster convergence in practice.
Specifically, if $F_L$ is nonexpansive, applying
Proposition~\ref{prop:km-intro} gives
\[ \left\| F(x_T)-x_T \right\|_2 = L\left\| F_L(x_T)-x_T \right\|_2\leqslant \frac{2L\left\| x_1-x_* \right\|_2}{\sqrt{T}},\quad T\geqslant 1. \]
for every fixed point $x_*$ of $F$. Focusing on the above bound, because
it scales with $L$, the best such guarantee would be obtained by the smallest
value of $L$ such that $F_L$ is nonexpansive.
This is closely related to the problem of
finding the right step-size (also called the learning rate, which here corresponds to
the inverse of $L$) in first-order optimization---the following
intuition is well-known to all practitioners, in particular in machine/deep learning: if the step-size is
too small, convergence happens but may be extremely slow, and if the
step-size is too large, convergence does not happen.

The above discussion can be extended to e.g.,\ positive definite
matrices $A\in \mathcal{S}_{++}^d(\mathbb{R})$ which can act as \emph{preconditioners}: operators
$F$ and $F_A:=I+A^{-1}(F-I)$ have the same fixed points.

Consider the following toy example. Let $\alpha>3$, $\varepsilon\in (0,1)$ and $F:\mathbb{R}^2\to \mathbb{R}^2$ be the linear
operator defined as
\[ F(x)=Mx, \quad x\in \mathbb{R}^2,\quad \text{where }M=\begin{pmatrix}-\alpha&0\\0&1-\varepsilon
\end{pmatrix}. \]
Operator $F$ has the origin as its unique fixed point because it is
linear and does not have $1$ as eigenvalue, and
because $-\alpha$ is the largest eigenvalue (in magnitude), $F$ is
$\alpha$-Lipschitz continuous for $\left\|\,\cdot\,\right\|_2$, with $\alpha$ being the optimal coefficient.
Because $\alpha>1$, operator $F$ is not nonexpansive, and the guarantee for
KM iterates from Proposition~\ref{prop:km-intro} does not apply. In fact,
starting from an initial point with nonzero first component,
both the power iterations $x_{t+1}=F(x_t)$ $(t\geqslant 1)$ and the standard KM
iterations $x_{t+1}=(x_t+F(x_t))/2$ can be seen to \emph{diverge}.

In this example however, it is easy to deduce that the operator 
\[ F_{(1+\alpha)/2}=I+\frac{2}{1+\alpha}\left( F-I \right),  \]
which corresponds to applying scaling $L=(\alpha+1)/2$ to $F$, 
\emph{is} nonexpansive. The corresponding KM iterations:
\[ x_{t+1}=\frac{x_t+F_{(1+\alpha)/2}(x_t)}{2},\quad t\geqslant 1, \]
enjoy the following convergence guarantee given by Proposition~\ref{prop:km-intro}:
\begin{equation}
\label{eq:3}
\left\| F(x_T)-x_T \right\|_2 =\frac{\alpha+1}{2}\left\| F_{(1+\alpha)/2}(x_T)-x_T \right\|_2 \leqslant \frac{(\alpha+1)\left\| x_1 \right\|_2}{\sqrt{T}},\quad T\geqslant 1.
\end{equation}

Furthermore, a better
convergence guarantee can be obtained by considering the following positive definite matrix:
\[ A=\begin{pmatrix}(1+\alpha)/2&0\\0&\varepsilon/2\end{pmatrix}. \]
The associated operator $F_A=I+A^{-1}(F-I)$
is also nonexpansive (because it is in fact equal to $-I$) and the corresponding KM iterations
$x_{t+1}^{(A)}=(x_t^{(A)}+F_A(x_t^{(A)}))/2$ ($t\geqslant 1$), with $x_1^{(A)}=x_1$,
satisfy the following convergence guarantee given again by Proposition~\ref{prop:km-intro}:
\[ \left\| A^{-1}(F(x_T^{(A)})-x_T^{(A)}) \right\|_2\leqslant\frac{2\left\| x_1 \right\|_2}{\sqrt{T}},\quad T\geqslant 1. \]
The above is a stronger guarantee than~\eqref{eq:3} because we can
deduce a bound for the second component of the residual that is
improved by a factor $(\alpha+1)/\varepsilon$. In this particular example, because
$F_A$ is simply $-I$, the above last guarantee is very loose, as the
iterates $(x_t^{(A)})_{t\geqslant 1}$ reach the fixed point in only one
iteration from any initial point, so the actual improvement is even stronger.

In general however, even when operator $F$ is given in an explicit form, it
may be very difficult to guess a positive scaling $L>0$ (resp.\ a positive definite
preconditioner $A$) such that $I+\frac{1}{L}(F-I)$ (resp. $I+A^{-1}(F-I)$)
is nonexpansive
and to benefit from the corresponding convergence guarantee. Even if initial operator $F$ is
already nonexpansive, there may be a scaling (resp.\ positive definite
preconditioner) that could provide a better convergence
guarantee, and that we might want to benefit from.
We therefore wish to define \emph{adaptive} iterations that enjoy
convergence guarantees that are similar to those given by the best
positive scaling (resp.\ the best positive definite preconditioner),
\emph{without prior knowledge} of the latter.

In the context of solving linear systems and optimization, a problem is said to be
\emph{ill-conditioned} when basic methods converge slowly because it
happens to be described in an unfavorable system of coordinates
\citep{nemirovsky1983problem,polyak1987introduction}. This issue is of
central importance in practice. Using available information on the
problem structure and geometry for finding, before performing the
iterations, a change of coordinates that hopefully accelerates
convergence, is called \emph{preconditioning}
\citep{young1954iterative,varga1962matrix,meijerink1977iterative}. The
approach we propose is to seek generic fixed-point iterations
that adapt \emph{on-the-fly} to the best-in-hindsight preconditioner
(in certain classes). This approach may be called \emph{online
  preconditioning},
 and has been very successful in first-order optimization through quasi-Newton
methods~\citep{broyden1970convergence,goldfarb1970family,fletcher1970new,shanno1970conditioning,davidon1991variable,fletcher1987practical}.

\subsection{Adversarial regret minimization}
\label{sec:regret-minimization}

This work introduces a new approach for defining and analyzing fixed-point iterations, based on adversarial regret minimization.
The latter is a class of sequential decision
problems where we aim at defining decision rules that provide worst-case guarantees
(i.e.,\ that hold for any sequence of decisions made by the
environment, seen as an adversary). The \emph{regret} is a
quantity that compares the actual gain obtained by the decision maker
with the gain of the best constant decision in hindsight (meaning the best gain
that \emph{could} have been obtained with a constant decision over
time with \emph{prior knowledge} of the decisions of the environment).

Adversarial regret minimization was first introduced in the fifties by
\citet{hannan1957approximation}, but underwent a spectacular
development over the last thirty years
\citep{cesa2006prediction,shalev2011online,hazan2016introduction,orabona2019modern},
with applications to numerous related
problems~\citep{freund1999adaptive,hart2000simple,cesa2004generalization,zinkevich2007regret,perchet2014approachability}.
In particular, deep links with first-order methods in convex and
nonconvex optimization were noticed and developed
\citep{levy2017online,cutkosky2019anytime,cutkosky2023optimal,defazio2024road,jiang2024improved}.
At the intersection with first-order optimization, a remarkable
breakthrough was the introduction of the AdaGrad family of algorithms
\citep{mcmahan2010adaptive,duchi2011adaptive}, which at their core are
 regret minimizers with theoretical guarantees that have
a strong \emph{adaptive} character, but then were found to also have
excellent performance in practice for first-order optimization, and in
particular, for the optimization of neural networks. As a matter of
fact, AdaGrad variants such as RMSprop~\citep{tieleman2012lecture} and
Adam~\citep{kingma2015adam} (and many others) are among the state-of-the-art
optimizers for deep learning.

\subsection{Contributions}
\label{sec:contribution}

The core contribution of this work is a conversion scheme of regret
minimization algorithms into fixed-point iterations, with convergence
guarantees derived from the regret bounds. This
connection builds upon the well-known equivalence between fixed points of
nonexpansive operators and zeros of co-coercive operators.
A notable special case is
the KM iterations which can be seen as the conversion of Online
Gradient Descent (one of the most basic regret minimization
algorithms) into fixed-point iterations, with the corresponding regret bound
recovering the optimal nonasymptotic bound (up to constant) for KM in
the Euclidean setting.\@

We then apply this conversion scheme to AdaGrad-type algorithms which
produce fixed-point iterations that enjoy adaptive convergence
guarantees of a new kind. Iterations based on AdaGrad-Norm (the
simplest version of AdaGrad, that involves a single data-dependent
step-size) are adaptive to the best-in-hindsight positive scaling,
meaning the lowest coefficient $L>0$ such that $I+\frac{1}{L}(F-I)$ is
nonexpansive. Going further, we introduce the notion of
$A$-nonexpansiveness (where $A$ is a positive definite matrix), which
corresponds to nonexpansiveness after applying preconditioner $A$, and
establish that iterations based on AdaGrad-Diagonal (resp.
AdaGrad-Full) are adaptive to the best-in-hindsight positive diagonal
(resp.\ positive definite) preconditioner. These can be seen as generic
variable metric methods that offer some on-the-fly remedy to
ill-conditioning. A few drawbacks should however be mentioned: (i)
those adaptive guarantees are on the best-iterate instead of the
last-iterate as in KM, (ii) dependence on the distance to the solution
is degraded unless an estimate is known, and (iii) in the case of
AdaGrad-Diagonal and AdaGrad-Full, the residual bounds depend on the
maximal distance of the iterates to a solution; a uniform rate follows
when the iterates are bounded, for example when the domain
is bounded.

\subsection{Related work}
\label{sec:related-work}

\subsubsection{Adaptive character of AdaGrad}
The AdaGrad family has been extensively studied in the context of
convex and nonconvex (possibly stochastic) first-order optimization
\citep{cutkosky2018distributed,ward2020adagrad,defossez2022simple,traore2021sequential,xie2020linear,faw2022power,liu2023convergence,attia2023sgd,wang2023convergence}.
Notably, AdaGrad-Norm was shown to be adaptive to the smoothness level
of the objective function (meaning the Lipschitz constant of the
gradient) in the convex case \citep{levy2018online}. The adaptivity to
positive scalings that we obtain for fixed points is analogous. In a
similar vein, AdaGrad-Diagonal has been proved in
\citet{liu2024large,jiang2024convergence} to be adaptive to the
anisotropic smoothness (aka coordinate-wise smoothness) of the
objective function. The adaptivity to the best-in-hindsight positive
diagonal preconditioner that we establish has a similar
spirit.

\subsubsection{Adaptive fixed-point iterations}
Various types of adaptivity in fixed-point iterations, different from
ours, have been proposed in the literature.
\citet{giselsson2016line} choose the step-size along $F(x_t)-x_t$
by testing residual decrease at candidate points.
Our AdaGrad iterations instead choose step-sizes from accumulated
residuals, without a line search.
SuperMann~\citep{themelis2019supermann} combines steps along proposed
directions, such as quasi-Newton directions, with acceptance tests
and fallback projections onto half-spaces containing the fixed points.
It guarantees global convergence, and superlinear convergence under
additional assumptions.
These guarantees differ from our regret-based comparison with
the best-in-hindsight scaling or preconditioner.
The closest related guarantee is provided by Algorithm 1 and Theorem 6 in
\citep{diakonikolas2020halpern}. These results use an initial scaling guess
$L_0>0$ and adapt to higher values only in an online fashion through a
doubling trick. They have the advantages over our guarantee in
Theorem~\ref{thm:adagrad-norm-fp} below of having the optimal dependence in
the distance to the solution and of the Halpern-like accelerated rate.
Theorem 4 in \citet{malitsky2020golden} obtains convergence to a fixed
point for demi-contractive operators, by estimating local Lipschitz
information online using the previous two iterates. Recently,
\citet{diakonikolas2025pushing} proposed Halpern-based iterations that
adapt to contraction, nonexpansiveness and even mild expansion.

\subsubsection{Monotone variational inequalities and games}
The present work builds on the equivalence between fixed points of
nonexpansive operators and zeros of co-coercive operators. Co-coercivity is
stronger than monotonicity, which brings to mind the topic of
variational inequalities with monotone operators. The latter are
typically solved using methods like
Extragradient~\citep{korpelevich1976extragradient} and its extension
Mirror-Prox~\citep{nemirovski2004prox} and have been related to regret
minimization~\citep{rakhlin2013optimization,chen2023generalized,sorin2024no,hadiji2024tracking},
in particular in the context of solving classes of games
~\citep{freund1999adaptive,mertikopoulos2019learning}. Adaptive
variants of such methods have also been proposed
~\citep{bach19a,diakonikolas2020halpern,antonakopoulos2021adaptive,ene2022adaptive,liu2022adaptive}.
Besides, extragradient methods have been related to a relaxation of
co-coercivity~\citep{gorbunov2022extragradient}. Thus, although regret
is known to be connected to various adjacent topics to nonexpansive
fixed-point iterations, it seems that the direct connection between regret
and fixed-point residuals of nonexpansive operators that we present in
this work, and its use to obtain adaptive fixed-point iterations,
have not been proposed before.

\subsubsection{Variable-metric methods}
Variable-metric fixed-point and splitting methods that are studied in
\citet{combettes2013variable,combettes2014variable} impose conditions
on a prescribed metric sequence. By contrast, the proposed
AdaGrad-type fixed-point iterations admit guarantees competing with
the best scaling, or diagonal/full-matrix preconditioning in hindsight.







\subsection{Disclosure of LLM use}
\label{sec:disclosure-llm-use}
The ideas leading to this work were developed over many years
starting in 2019. The first preprint appeared online in September
2025. LLMs were only used for proofreading and polishing the revised
version of the paper.

\subsection{Notation}
Throughout the paper, \(d\geqslant 1\) is an integer, \(\mathcal{X}\subset \mathbb{R}^d\) is a 
closed convex set that is neither empty nor a singleton.
\(I\) denotes the identity operator whose domain is to be deduced from the
context (usually \(\mathcal{X}\) or \(\mathbb{R}^d\)).
The canonical inner product and Euclidean norm in $\mathbb{R}^d$ are denoted as
\[ \left< x, x'\right> =\sum_{i=1}^dx_ix_i'\quad\text{and}\quad \left\| x \right\|_2 =\sqrt{\left< x, x \right> },\quad x,x'\in \mathbb{R}^d. \]
\(I_d\) denotes the identity matrix of size \(d\).
For a vector $u\in \mathbb{R}^d$, $\operatorname{diag}u$ denotes the 
diagonal matrix of size $d$ with diagonal coefficients given by $u$.
If $A$ is a symmetric positive semidefinite matrix, its square root is denoted $A^{1/2}$.
The set of symmetric positive definite matrices of size $d$ is denoted $\mathcal{S}_{++}^d(\mathbb{R})$.
For a symmetric positive semidefinite matrix \(A\), define the seminorm
\[ \left\| x \right\|_A=\sqrt{\left< x, Ax \right>}. \]
For \(A\in \mathcal{S}_{++}^d(\mathbb{R})\), the associated dot product and Mahalanobis norm are denoted:
\[ \left< x, x' \right>_A=\left< x, Ax' \right> \quad \text{and}\quad  \left\| x \right\|_A=\sqrt{\left< x, x \right>_A},\quad x,x'\in \mathbb{R}^d. \]
For a given norm $\left\|\,\cdot\,\right\|$ in $\mathbb{R}^d$, the (possibly
infinite) diameter of a set $\mathcal{A}\subset \mathbb{R}^d$ is defined as
\[ \operatorname{diam}_{\left\|\,\cdot\,\right\| }(\mathcal{A})=\sup_{x,x'\in \mathcal{A}}\left\| x'-x \right\|. \]
The diameter with respect to $\left\|\,\cdot\,\right\|_2$ and
$\left\|\,\cdot\,\right\|_{\infty}$ are denoted $\operatorname{diam}_2$ and
$\operatorname{diam}_{\infty}$ respectively.
The Euclidean projection (resp.\ the projection with respect to
Mahalanobis norm $\left\|\,\cdot\,\right\|_A$ for a given matrix
$A\in \mathcal{S}_{++}^d(\mathbb{R})$) onto a closed convex set $\mathcal{X}\subset \mathbb{R}^d$ is denoted $\Pi_{\mathcal{X}}$
(resp.\ $\Pi_{\mathcal{X},A}$). The unit simplex in $\mathbb{R}^d$ is denoted by
\[ \Delta_d=\left\{ x\in \mathbb{R}_+^d,\ \sum_{i=1}^dx_i=1 \right\}.  \]

\section{Properties of operators}
\label{sec:basic-prop-oper}

This section recalls several basic definitions and introduces
$A$-nonexpansiveness. A key remark is that operators $F$ and
$F_L:=I+\frac{1}{L}(F-I)$ (resp.\ $F_A:=I+A^{-1}(F-I)$)
have the same fixed points
for any $L>0$ (resp.\ for any matrix
$A\in \mathcal{S}_{++}^d(\mathbb{R})$). Note that $F_L$ is a special case of $F_A$ (with $A=LI_d$).
Operator $F_L$ can be interpreted as a scaled version of $F$, and 
operator $F_A$ is obtained from $F$ via a kind of preconditioning.

Even if $F$ is not nonexpansive,
$F_A$ may very well be nonexpansive (for e.g.,\ some Mahalanobis norm) and KM
iterations with $F_A$ would then converge to a fixed point of $F$.
Even if $F$ \emph{is} nonexpansive, $F_A$ may yield a better
convergence for some matrix $A$. It will turn out that a natural
notion is $F_A$ being nonexpansive with respect to
$\left\|\,\cdot\,\right\|_A$, and that will be the definition of
$A$-nonexpansiveness. Besides, the fixed points of an operator $F$ are
the zeros of associated operator $G=\frac{I-F}{2}$, and the
characterization through $G$ of ($A$-)nonexpansiveness of $F$ will be
particularly useful.
\begin{definition}
  Let \(F,G:\mathcal{X}\to \mathbb{R}^d\) and \(x_*\in \mathcal{X}\).
\begin{enumerate}[(i)]
\item \(x_*\in \mathcal{X}\) is a \emph{fixed point} of operator \(F\) if \(F(x_*)=x_*\).
\item \(x_*\in \mathcal{X}\) is a \emph{zero} of operator \(G\) if \(G(x_*)=0\).
\end{enumerate}
\end{definition}

We now recall nonexpansiveness and introduce two successive extensions involving a
positive scalar and a positive definite matrix respectively.
\begin{definition}[Nonexpansiveness]
Let  \(F:\mathcal{X}\to \mathbb{R}^d\) and \(\left\|\,\cdot\,\right\| \) a norm in
\(\mathbb{R}^d\).
\begin{enumerate}[(i)]
\item Operator \(F\) is \emph{nonexpansive} with respect to \(\left\|\,\cdot\,\right\|\) if for all
\(x,x'\in \mathcal{X}\),
\[ \left\| F(x')-F(x) \right\| \leqslant \left\| x'-x \right\|. \]
If $F$ is nonexpansive with respect to $\left\|\,\cdot\,\right\|_2$, $F$ is
simply said to be nonexpansive.
\item Let $L>0$. Operator $F$ is \emph{$L$-nonexpansive} if
  $I+\frac{1}{L}(F-I)$ is nonexpansive for $\left\|\,\cdot\,\right\|_2$.
\item Let $A\in \mathcal{S}_{++}^d(\mathbb{R})$. Operator $F$ is \emph{$A$-nonexpansive} if
  $I+A^{-1}(F-I)$ is nonexpansive for $\left\|\,\cdot\,\right\|_A$.
\end{enumerate}
\end{definition}
In the case where $L\in (0,1)$, $L$-nonexpansiveness corresponds to
being \emph{$L$-averaged}, which is a notion that was first formalized
in \citet{baillon1978asymptotic}. It is called so because the operator
then is a convex combination of the identity and a nonexpansive
operator. We here consider any value $L>0$.

Obviously, the notion of $A$-nonexpansiveness encompasses the other
two, which are recovered by $A=I_d$ and $A=LI_d$ respectively.
Note that $A$-nonexpansiveness of $F$ is \emph{not} equivalent to nonexpansiveness of
$F$ with respect to $\left\|\,\cdot\,\right\|_A$. As noted in Section~\ref{sec:introduction}, an operator $F$ may fail to be nonexpansive and yet be
$A$-nonexpansive for some matrix $A\in \mathcal{S}_{++}^d(\mathbb{R})$.

We now recall the notion of co-coercivity, and allow for arbitrary norms. It is well-known that the
nonexpansiveness of $F$ (with respect to $\left\|\,\cdot\,\right\|_2$) is characterized by the co-coercivity of
$G=\frac{I-F}{2}$ (with respect to $\left\|\,\cdot\,\right\|_2$). Regarding $A$-nonexpansiveness, its
characterization will naturally involve co-coercivity with respect to $\left\|\,\cdot\,\right\|_A$.
\begin{definition}[Co-coercivity for arbitrary norms]
Let \(G:\mathcal{X}\to \mathbb{R}^d\), \(\left\|\,\cdot\,\right\|\) a norm in \(\mathbb{R}^d\), and \(L>0\).
  Operator \(G\) is \(L\)-\emph{co-coercive} with respect to \(\left\|\,\cdot\,\right\|\) if for all \(x,x'\in \mathcal{X}\),
  \[ \left< G(x')-G(x), x'-x \right> \geqslant \frac{1}{L}\left\| G(x')-G(x) \right\|^2_*, \]
  where $\left\|\,\cdot\,\right\|_*$ denotes the dual norm of $\left\|\,\cdot\,\right\|$.
\end{definition}
In the literature, \(L\)-co-coercivity is sometimes defined with
factor \(L\) in the above right-hand side instead of \(1/L\)---see
e.g.,\ \citep{bauschke2011convex}.
The above convention has the advantage of the following linear property: if
\(G\) is \(L\)-co-coercive and \(\eta>0\), then \(\eta G\) is \(\eta L\)-co-coercive.

Let us keep in mind that the same solutions can be expressed through
various operators of interest: for \(F:\mathcal{X}\to \mathbb{R}^d\) and
\(A\in \mathcal{S}_{++}^d(\mathbb{R})\), a point
\(x_*\in \mathcal{X}\) is a fixed point of \(F\) if, and only if, it is a fixed
point of \(I+A^{-1}(F-I)\), which is also equivalent to being a zero of
$G:=\frac{I-F}{2}$.

We now state a natural extension of the well-known characterization of nonexpansiveness (see
e.g.,\ \citep[Proposition 4.2]{bauschke2011convex}) to $A$-nonexpansiveness. 
\begin{proposition}
  \label{prop:characterization}
  Let \(F:\mathcal{X}\to \mathbb{R}^d\), \(G=\frac{I-F}{2}\), and
  \(A\in \mathcal{S}_{++}^d(\mathbb{R})\).
  Then, $F$ is
  $A$-nonexpansive if, and only if, \(G\) is $1$-co-coercive for
  $\left\|\,\cdot\,\right\|_A$, in other words,
\begin{equation}
\label{eq:1}
\forall x,x'\in \mathcal{X},\quad \left< G(x')-G(x), x'-x \right>\geqslant \left\| G(x')-G(x) \right\|^2_{A^{-1}}. 
\end{equation}
\end{proposition}
\begin{proof}
  Note that \(I+A^{-1}(F-I)=I-2A^{-1}G\). This operator being
  nonexpansive for \(\left\|\,\cdot\,\right\|_A\) writes
  \[ \left\| (I-2A^{-1}G)(x')-(I-2A^{-1}G)(x) \right\|_A^2\leqslant \left\| x'-x \right\|_A^2,\quad x,x'\in \mathcal{X}. \]
Let \(x,x'\in \mathcal{X}\) be given. Expanding the above left-hand side, we obtain
\begin{multline*}
\left\| (I-2A^{-1}G)(x')-(I-2A^{-1}G)(x) \right\|_A^2=4\left\| A^{-1}(G(x')-G(x)) \right\|_A^2\\
\qquad -4\left< A^{-1}(G(x')-G(x)), A(x'-x) \right> + \left\| x'-x \right\|_A^2\\
=4\left\| G(x')-G(x) \right\|^2_{A^{-1}}-4\left< G(x')-G(x), x'-x \right>+\left\| x'-x \right\|_A^2. 
\end{multline*}
Plugging the above into the first inequality and rearranging gives the result.
\end{proof}
Through the lens of the above characterization, the notion
$A$-nonexpansiveness may feel more natural: for operator $F$, going
from nonexpansiveness to $A$-nonexpansiveness corresponds to norm
$\left\|\,\cdot\,\right\|_2$ being simply replaced by
$\left\|\,\cdot\,\right\|_A$ in the co-coercivity of associated operator
$G=\frac{I-F}{2}$. Yet another equivalent formulation is $A^{-1}G$
being firmly nonexpansive for $\left\|\,\cdot\,\right\|_A$.

\section{Fixed points via Online Linear Optimization}
\label{sec:regr-minim-key}

This section first recalls the classical adversarial regret
minimization framework called Online Linear Optimization as well as the
 regret bound of Online Gradient Descent (OGD), which is one of the
simplest algorithms. We then turn to the connection between fixed
points and regret: Lemma~\ref{lm:online-to-fp} provides an upper bound
on fixed-point residuals that can be interpreted as a regret. This
connection naturally suggests a scheme for converting regret
minimization algorithms into fixed-point iterations, because bounds on
the regret would then become bounds on the fixed-point residuals. As
an illustration, Corollary~\ref{cor:km} applies this conversion
scheme to Online Gradient Descent and recovers the fundamental
KM iteration and its nonasymptotic analysis in the Euclidean case.
Then, considering more general nonexpansive non-self mappings
$F:\mathcal{X}\to \mathbb{R}^d$, we study a projection variant of KM, obtained by converting
Projected Online Gradient Descent.

\subsection{Online linear optimization}
\label{sec:online-line-optim-1}
Online Linear Optimization \citep{zinkevich2003online,kalai2005efficient} can be described as a sequential
decision problem where the \emph{Decision Maker} chooses its \emph{actions} in the
closed convex set $\mathcal{X}\subset \mathbb{R}^d$: at step $t\geqslant 1$,
\begin{itemize}
\item the Decision Maker chooses \emph{action} $x_t\in \mathcal{X}$,
\item Nature chooses and reveals \emph{payoff vector} $u_t\in \mathbb{R}^d$,
\item the Decision Maker obtains \emph{payoff} $\left< u_t, x_t \right>$.
\end{itemize}
Let $T\geqslant 1$ be an integer. The Decision Maker wishes to maximize its
cumulative payoff $\sum_{t=1}^T\left< u_t, x_t \right>$, but no worst-case type guarantee can
be achieved on this quantity, in an absolute sense. We instead
aim at establishing upper bounds on the \emph{regret} (with respect to some comparison point $x\in \mathcal{X}$):
\begin{equation}
\label{eq:regret}
\sum_{t=1}^T\left< u_t, x \right> - \sum_{t=1}^T\left< u_t, x_t \right>=\sum_{t=1}^T\left< u_t, x-x_t \right>,  
\end{equation}
which compares the cumulative payoff actually obtained by the Decision
Maker to the cumulative payoff that it would have obtained by choosing
action $x$ at each step. We recall the elementary regret bound guaranteed by
OGD in the unconstrained case $\mathcal{X}=\mathbb{R}^d$.
\begin{proposition}[Unconstrained Online Gradient Descent]
  \label{prop:ogd}
 Let \((u_t)_{t\geqslant 1}\) be a sequence in \(\mathbb{R}^d\), $\eta>0$, \(x_1\in \mathbb{R}^d\) and
 \begin{equation}
\label{eq:ogd}
x_{t+1}=x_t+\eta u_t,\quad t\geqslant 1. 
\end{equation}
Then for every integer \(T \geqslant 1\) and $x\in \mathbb{R}^d$,
 \[ \sum_{t=1}^{T}\left< u_t, x-x_t \right> \leqslant \frac{\left\| x-x_1 \right\|_2^2}{2\eta}+\frac{\eta}{2}\sum_{t=1}^{T}\left\| u_t \right\|_2^2. \]
\end{proposition}
\begin{proof}
For $t\geqslant 1$, writing
\[ \left\| x_{t+1}-x \right\|_2^2=\left\| x_t+\eta u_t-x \right\|_2^2 =\left\| x_t-x \right\|_2^2+2\left< \eta u_t, x_t-x \right> + \eta^2\left\| u_t \right\|_2^2,\]
then dividing by $2\eta$, summing over $t=1,\dots,T$ and rearranging
gives the result.
\end{proof}
If the vectors $(u_t)_{t\geqslant 1}$ are assumed to be
bounded, choosing $\eta=1/\sqrt{T}$ above yields a regret bound that grows as $\sqrt{T}$.
In other words, the \emph{average} regret has an upper bound that
vanishes as $1/\sqrt{T}$, which is a foundational result in
adversarial regret minimization: with such an algorithm, for large $T$ and on average, the
Decision Maker is guaranteed to perform as well as the constant
algorithm that chooses action $x$ at each step, this being true for all $x\in \mathcal{X}$.

Beyond sequential decision problems \emph{per se}, regret minimization
has deep connections with first-order optimization. OGD is of course
an extension of Gradient Descent. But many other optimization
algorithms can be analyzed using regret bounds e.g.,\ the Nesterov
Accelerated Gradient method~\citep{nesterov1983method} and Mirror
Descent~\citep{nemirovsky1983problem,beck2003mirror}. We now turn to
the connection between regret and fixed points.

\subsection{Key lemma}
\label{sec:regr-nonexp-oper}

The following lemma is the key connection between fixed points and
regret for $A$-nonexpansive operators.

It provides an upper bound on a squared \emph{fixed-point residual}, which is
a measure of how far a point is from being a fixed point---in
particular, the residual is zero at a point if, and only if, it is a
fixed point. We then remark that this upper bound can be interpreted
as a (one-step) regret.
\begin{lemma}
  \label{lm:online-to-fp}
  Let $A\in \mathcal{S}_{++}^d(\mathbb{R})$,
  \(F:\mathcal{X}\to \mathbb{R}^d\) an $A$-nonexpansive operator and
  \(x_*\in \mathcal{X}\) a fixed point of \(F\). Then for all \(x\in \mathcal{X}\),
  \[ \left\| F(x)-x \right\|_{A^{-1}}^2\leqslant 2\left< F(x)-x, x_*-x \right>. \]
\end{lemma}
\begin{proof}
  Let \(G=\frac{I-F}{2}\), which according to Proposition~\ref{prop:characterization} is
  \(1\)-co-coercive for \(\left\|\,\cdot\,\right\|_A\).
  By definition of \(x_*\), \(G(x_*)=0\).
  Let \(x\in \mathcal{X}\). Using co-coercivity,
  \begin{align*}
  \left\| G(x) \right\|_{A^{-1}}^2&=\left\| G(x)-G(x_*) \right\|_{A^{-1}}^2\leqslant \left< G(x)-G(x_*), x-x_* \right> \\
    &=\left< G(x), x-x_* \right>=\left< \frac{x-F(x)}{2}, x-x_* \right>. 
  \end{align*}
  Besides, \(\left\| G(x) \right\|_{A^{-1}}^2\) rewrites as
  \[ \left\| G(x) \right\|_{A^{-1}}^2=\left\| \frac{x-F(x)}{2} \right\|_{A^{-1}}^2=\frac{1}{4}\left\| F(x)-x \right\|_{A^{-1}}^2. \]
  The result follows from simplifying.
\end{proof}
We deduce the following immediate corollary by considering the
weighted sum of the inequality from the above lemma, for each of the
three cases: nonexpansiveness, $L$-nonexpansiveness and $A$-nonexpansiveness.
\begin{corollary}
  \label{cor:online-to-fp}
  Let \(F:\mathcal{X}\to \mathbb{R}^d\), \(x_*\in \mathcal{X}\) a fixed point of \(F\), 
  \(T\geqslant 1\) an integer, \(x_1,\dots,x_T\in \mathcal{X}\), and \(\gamma_1,\dots,\gamma_T>0\).
\begin{enumerate}[(i)]
\item\label{item:online-to-fp-nonexp} If \(F\) is nonexpansive, then
\[ \sum_{t=1}^T\gamma_t\left\| F(x_t)-x_t \right\|_2^2\leqslant 2\sum_{t=1}^T\left< \gamma_t(F(x_t)-x_t), x_*-x_t \right>.  \]
\item\label{item:online-to-fp-nonexp-L} Let \(L>0\). If $F$ is
  $L$-nonexpansive, then
\[ \sum_{t=1}^T\gamma_t\left\| F(x_t)-x_t \right\|_2^2\leqslant 2L\sum_{t=1}^T\left< \gamma_t(F(x_t)-x_t), x_*-x_t \right>.  \]
\item\label{item:online-to-fp-nonexp-A} Let \(A\in \mathcal{S}_{++}^d(\mathbb{R})\). If
  $F$ is $A$-nonexpansive, then
\[ \sum_{t=1}^T\gamma_t\left\| F(x_t)-x_t \right\|^2_{A^{-1}}\leqslant 2\sum_{t=1}^T\left< \gamma_t(F(x_t)-x_t), x_*-x_t \right>.  \]
\end{enumerate}
\end{corollary}
By considering payoff vectors $u_t=\gamma_t(F(x_t)-x_t)$ for $t\geqslant 1$, the above upper
bounds have the form of a cumulative regret as in~\eqref{eq:regret} and
thus suggest the following approach. If points $(x_t)_{t\geqslant 1}$ are
chosen so that some upper bound on the regret is guaranteed, the same
bound holds for the (weighted) sum of the squared fixed-point residuals.
Interestingly, in cases~\eqref{item:online-to-fp-nonexp-L} and~\eqref{item:online-to-fp-nonexp-A}, the expression of the cumulative regret
$\sum_{t=1}^T\left< \gamma_t(F(x_t)-x_t), x_*-x_t \right>$ does not depend on
$L$ or $A$, and therefore, the regret minimizing
algorithm that chooses the points $(x_t)_{t\geqslant 1}$ \emph{need not have prior
knowledge} of $L$ or $A$. This feature will be of importance in
Sections~\ref{sec:adapt-co-coerc} and \ref{sec:adapt-posit-defin} when obtaining
guarantees that adapt to $L$-/$A$-nonexpansiveness \emph{without prior knowledge} of $L$ or $A$.

The nonexpansiveness assumptions in Lemma~\ref{lm:online-to-fp} and
Corollary~\ref{cor:online-to-fp} above can easily be relaxed into \emph{quasi-nonexpansiveness}.
For a given fixed point $x_*\in \mathcal{X}$ of $F$ and $A\in \mathcal{S}_{++}^d(\mathbb{R})$, we can
define $A$-quasi-nonexpansiveness of operator $F$ with respect to
$x_*$ as operator $F_A:=I+A^{-1}(F-I)$ not increasing the distance to
$x_*$ (measured with $\left\|\,\cdot\,\right\|_A$), in other words
\[ \forall x\in \mathcal{X},\quad \left\| F_A(x)-x_* \right\|_A \leqslant \left\| x-x_* \right\|_A. \]
Then, extending Proposition~\ref{prop:characterization}, this property
can be seen to be equivalent to \emph{star-co-coercivity} of $G$ (see~\citep{gorbunov2022extragradient}):
\[\forall x\in \mathcal{X},\quad  \left\| G(x) \right\|_{A^{-1}}^2\leqslant \left< G(x), x-x_* \right>, \]
where $G=(I-F)/2$, and the conclusions of Lemma~\ref{lm:online-to-fp} and
Corollary~\ref{cor:online-to-fp} then remain true. The details are worked out in
Appendix~\ref{sec:extens-quasi-nonexp}.

Regarding $L$-(quasi-)nonexpansiveness, we can further relax the
\emph{global} character of the assumption, by considering in the lemma
below a \emph{local} coefficient \(L_T\) that only depends on points
\(x_*,x_1,\dots,x_T\) and their images by operator \(F\).
\begin{lemma}
  \label{lm:online-to-fp-local}
  Let  \(F:\mathcal{X}\to \mathbb{R}^d\) be an operator, 
  \(x_*\in \mathcal{X}\), \(T\geqslant 1\) an integer and \(x_1,\dots,x_T\in \mathcal{X}\) such that
  \[ L_T:=\inf_{}\left\{ L>0,\ \forall 1\leqslant t\leqslant T,\ \left\| F(x_t)-x_t \right\|_2^2\leqslant 2L\left< F(x_t)-x_t, x_*-x_t \right> \right\}  \]
  is finite.
\begin{enumerate}[(i)]
\item\label{item:1} Then,
\[ \sum_{t=1}^T\left\| F(x_t)-x_t \right\|_2^2\leqslant 2L_T\sum_{t=1}^T\left< F(x_t)-x_t, x_*-x_t \right>.  \]
\item\label{item:2} Let $L>0$. If $F$ is $L$-nonexpansive and if $x_*$ is a fixed point
  of $F$, then $L_T\leqslant L$.
\end{enumerate}
\end{lemma}
\begin{proof}
  \eqref{item:1} follows from the definition of $L_T$ and
  \eqref{item:2} follows from Lemma~\ref{lm:online-to-fp}.
\end{proof}

If $x_*\in \mathcal{X}$ is a fixed point of $F$ and operator $F$ does satisfy
$L$-nonexpansiveness for some \emph{global} coefficient $L>0$, then the above
\emph{local} coefficient satisfies $L_T\leqslant L$ by
Proposition~\ref{prop:characterization} and we recover case~\eqref{item:online-to-fp-nonexp-L} in
Corollary~\ref{cor:online-to-fp}. But in cases where global coefficient $L$
is infinite, or where local coefficient $L_T$ is much lower, the above lemma will allow in
Section~\ref{sec:adapt-co-coerc} to
obtain guarantees that scale with $L_T$.

\subsection{Krasnoselskii--Mann iterations as a special case}
\label{sec:krasnoselskii-mann-as}

As a basic illustration of our conversion scheme, we combine
Lemma~\ref{lm:online-to-fp} with the Online Gradient Descent
algorithm~\eqref{eq:ogd} to reinterpret the classical KM iterations
and its nonasymptotic analysis in the Euclidean/Hilbert
setting~\citep{vaisman2005convergencia,cominetti2014rate}\footnote{\emph{Asymptotic}
  rates are faster---see
  \citet{browder1966solution,matsushita2017convergence}.}.
\begin{corollary}[Krasnoselskii--Mann iterations]
  \label{cor:km}
  Let \(F:\mathcal{X}\to \mathcal{X}\) be a nonexpansive operator, $x_*\in \mathcal{X}$ a fixed point of $F$,
  \((\gamma_t)_{t\geqslant 1}\) a sequence in \((0,1)\),
  \(x_1\in \mathcal{X}\) and for \(t\geqslant 1\),
  \[ x_{t+1}=\gamma_tF(x_t)+(1-\gamma_t)x_t. \]
  Then for all \(T\geqslant 1\),
\[ \left\| F(x_T)-x_T \right\|_2\leqslant \frac{\left\| x_1-x_* \right\|_2}{\sqrt{\sum_{t=1}^T(1-\gamma_t)\gamma_t}}. \]
\end{corollary}
\begin{proof}
The iteration rewrites
\[ x_{t+1}=x_t+\gamma_t(F(x_t)-x_t),\quad t\geqslant 1, \]
which corresponds to Online Gradient Descent~\eqref{eq:ogd} with
payoff vectors \(u_t=\gamma_t(F(x_t)-x_t)\) (for \(t\geqslant 1\)).
Combining the regret bound from Proposition~\ref{prop:ogd} with
bound~\eqref{item:online-to-fp-nonexp} in 
Corollary~\ref{cor:online-to-fp} gives
\[ \sum_{t=1}^T\gamma_t\left\| F(x_t)-x_t \right\|_2^2\leqslant \left\| x_1-x_* \right\|_2^2+\sum_{t=1}^T\gamma_t^2\left\| F(x_t)-x_t \right\|_2^2,  \]
which simplifies into
\[ \sum_{t=1}^T(1-\gamma_t)\gamma_t\left\| F(x_t)-x_t \right\|_2^2\leqslant \left\| x_1-x_* \right\|_2^2. \]
Let us now prove that \((\left\| F(x_t)-x_t \right\|_2)_{t\geqslant 1}\) is
nonincreasing, the result will follow. For \(t\geqslant 1\), using the
nonexpansiveness of \(F\),
\begin{align*}
\left\| F(x_{t+1})-x_{t+1} \right\|_2&=\left\| F(x_{t+1})-F(x_{t})+F(x_{t})-x_{t+1} \right\|_2 \\
  &\leqslant \left\| F(x_{t+1})-F(x_{t}) \right\|_2 + \left\| F(x_{t})-x_{t+1} \right\|_2 \\
  &\leqslant \left\| x_{t+1}-x_{t} \right\|_2 + (1-\gamma_t)\left\| F(x_{t})-x_{t} \right\|_2 \\
  &=\gamma_t\left\| F(x_{t})-x_{t} \right\|_2 + (1-\gamma_t)\left\| F(x_{t})-x_{t} \right\|_2\\
  &=\left\| F(x_{t})-x_{t} \right\|_2.
\end{align*}
\end{proof}

\subsection{Projected Krasnoselskii--Mann iterations}
\label{sec:proj-krasn-mann}

A requirement of KM is that
$F:\mathcal{X}\to \mathcal{X}$ is a \emph{self mapping}, meaning its images belong to
its convex domain $\mathcal{X}$, so that the next iterate can be chosen as a
weighted average of the current iterate and its image, which then
remains in the domain.

However, the discussion from Section~\ref{sec:quest-adaptive-fixed}
that compares the guarantee offered by KM applied to $F$ and
$F_L:=I+\frac{1}{L}(F-I)$ respectively does not always make sense
because for $L<1$, $F_L$ is not necessarily a self mapping and KM
iterations with $F_L$ may not be defined. Those properties are true
for instance if $\mathcal{X}=\mathbb{R}^d$. We now consider a \emph{projected} variant of
KM that is always defined.

For the remainder of the paper, we remove the self mapping requirement
and consider operators $F:\mathcal{X}\to \mathbb{R}^d$ that admit a fixed point
$x_*\in \mathcal{X}$. Then, KM iterations are not always defined, but
converting the \emph{Projected} Online Gradient
Descent~\citep{zinkevich2003online}, we obtain below a projected
variant of KM that enjoys the same residual bound as in
Corollary~\ref{cor:km}, albeit on the best-iterate, instead of the
last-iterate.\@

\begin{theorem}[Projected Krasnoselskii--Mann iterations]
\label{thm:proj-KM}
  Let \(F:\mathcal{X}\to \mathbb{R}^d\) be a nonexpansive operator, $x_*\in \mathcal{X}$ a fixed point
  of $F$,
  \((\gamma_t)_{t\geqslant 1}\) a sequence in \((0,1)\),
  \(x_1\in \mathcal{X}\) and for \(t\geqslant 1\),
  \[ x_{t+1}=\Pi_{\mathcal{X}}(\gamma_tF(x_t)+(1-\gamma_t)x_t). \]
  Then for all \(T\geqslant 1\),
\[ \min_{1\leqslant t\leqslant T}\left\| F(x_t)-x_t \right\|_2\leqslant \frac{\left\| x_1-x_* \right\|_2}{\sqrt{\sum_{t=1}^T(1-\gamma_t)\gamma_t}}. \]
\end{theorem}
\begin{proof}
For $t\geqslant 1$, considering payoff vector $u_t=\gamma_t(F(x_t)-x_t)$, the points
$(x_t)_{t\geqslant 1}$ are the corresponding Projected Online Gradient Descent
iterates, and we can write
\begin{align*}
\sum_{t=1}^T\gamma_t\left\| F(x_t)-x_t \right\|_2^2&\leqslant 2\sum_{t=1}^T\left< \gamma_t(F(x_t)-x_t), x_*- x_t \right>\\
  &\leqslant \left\| x_1-x_* \right\|_2^2+\sum_{t=1}^T\gamma_t^2\left\| F(x_t)-x_t\right\|_2^2,
\end{align*}
where the first inequality holds by Corollary~\ref{cor:online-to-fp}
and the second inequality is by the regret bound from
Proposition~\ref{prop:proj-ogd}. Rearranging gives
\[ \sum_{t=1}^T\gamma_t(1-\gamma_t)\left\| F(x_t)-x_t \right\|_2^2\leqslant \left\| x_1-x_* \right\|_2^2. \]
Dividing by $\sum_{t=1}^T\gamma_t(1-\gamma_t)$ yields
\[ \min_{1\leqslant t\leqslant T}\left\| F(x_t)-x_t \right\|_2^2\leqslant \frac{\sum_{t=1}^T\gamma_t(1-\gamma_t)\left\| F(x_t)-x_t \right\|_2^2}{\sum_{t=1}^T\gamma_t(1-\gamma_t)}\leqslant \frac{\left\| x_1-x_* \right\|_2^2}{\sum_{t=1}^T(1-\gamma_t)\gamma_t}, \]
hence the result.
\end{proof}

Like for basic KM iterations in Corollary~\ref{cor:km}, because the
upper bound is the same, the best guarantee is obtained by choosing
$\gamma_t=1/2$ for all $t\geqslant 1$. However, the bound holds for the
best iterate instead of the last iterate as in
Corollary~\ref{cor:km}. In the context of smooth convex optimization,
the corresponding algorithm is Projected Gradient Descent, which is
known to offer a last-iterate guarantee on the optimality gap
$f(x_T)-f(x_*)$. Whether an analogous last-iterate guarantee is
possible for fixed-point iterations is an interesting question.




There are other approaches for handling non-self mappings.
One idea is to consider operator $\Pi_{\mathcal{X}}\circ F$ which is a self-mapping,
is nonexpansive as soon as $F$ is (because the projection is also
nonexpansive), and has the same fixed points (provided that $F$ is
nonexpansive and has a fixed point in $\mathcal{X}$). However, the drawback of
this approach is that one naturally obtains guarantees on
$\left\| \Pi_{\mathcal{X}}(F(x))-x \right\|_2$ instead of the original residual
$\left\| F(x)-x \right\|_2$, and controlling the latter quantity with
the former involves an undesirable factor $\left\| x-x_* \right\|_2$. Alternatively,
\citet{colao2015krasnoselskii} proposed to choose the averaging
weights of the KM iteration at each step in such a way that the next
iterate \emph{does} belong to the domain $\mathcal{X}$---in addition to
nonexpansiveness of operator $F$, geometric conditions on the domain
$\mathcal{X}$ (strict convexity) and the operator $F$ (inward condition) are
required.

\section{AdaGrad-type iterations}
\label{sec:adagr-type-iter}
The AdaGrad family was initially introduced as regret
minimization algorithms with step-sizes that depend on data observed in
previous steps in a specific
way~\citep{mcmahan2010adaptive,duchi2011adaptive}. Perhaps the most
popular variant is called AdaGrad-Diagonal, and has per-coordinate
step-sizes. This variant has been very successful for first-order
optimization, especially in neural network optimization. Nowadays, the
state-of-the-art optimizers for deep learning are variants of
AdaGrad-Diagonal with stronger performance but poorer theoretical
understanding, e.g.,\ \citet{tieleman2012lecture,kingma2015adam}. Another
variant within the family is AdaGrad-Full, which involves scalings using
full matrices (instead of being restricted to diagonal matrices in the case of
per-coordinate step-sizes), but has been less popular in practice,
because each iteration then requires the computation of $A_t^{-1}$
multiplied by a vector, where $A_t$ is a full $d\times d$ matrix, which
does not scale to high dimensions.

We first focus on the simplest variant, sometimes called
AdaGrad-Norm, which has a single step-size that depends on the data
from previous steps. In the context of first-order convex
optimization, an important theoretical guarantee explaining the good
practical performance was obtained by \citet{levy2018online}. It established that
the algorithm was adaptive to the unknown Lipschitz constant of the
gradient of the objective function (also called the smoothness
constant), in other words, it obtained a
convergence guarantee that was similar to Gradient Descent tuned with
prior knowledge of the smoothness constant. The adaptivity to unknown
$L$-nonexpansiveness that we obtain for fixed points in
Theorem~\ref{thm:adagrad-norm-fp} is analogous to that result, but is not
a direct extension, as will be discussed.

AdaGrad-Diagonal and AdaGrad-Full will also be converted for fixed points in
Section~\ref{sec:adapt-posit-defin} and will enjoy guarantees with
adaptivity to unknown $A$-nonexpansiveness.

\subsection{Adaptivity to positive scalings}
\label{sec:adapt-co-coerc}


For $L>0$, operator $F$ being $L$-nonexpansive means
by definition that associated operator $F_L:=I+\frac{1}{L}(F-I)$ is
nonexpansive. If such a coefficient $L$ were known, and if operator
$F$ admits a fixed point $x_*\in \mathcal{X}$, KM (resp.\ Projected KM, if
$\mathcal{X}$ is a strict subset of $\mathbb{R}^d$) with operator
$F_L$ could be used (with e.g.,\ constant weights $\gamma_t=1/2$,
$t\geqslant 1$) and would offer the residual bound
$2L\left\| x_1-x_* \right\|_2/\sqrt{T}$ on the last iterate
(resp.\ best iterate). Theorem~\ref{thm:adagrad-norm-fp} below gives a
best-iterate \emph{adaptive} guarantee with the same dependence on
$L$ and $T$ (up to logarithms) without prior knowledge of either.
Besides, the
linear dependence on $\left\| x_1-x_* \right\|_2$ is maintained if a
good estimate of this distance is known, and degrades to a quadratic
dependence otherwise.

For a sequence \((u_t)_{t\geqslant 1}\) in \(\mathbb{R}^d\), parameter $\eta >0$ and initial point $x_1\in \mathcal{X}$,
the AdaGrad-Norm iterates on $\mathcal{X}$, first studied in \citet{levy2017online}, are defined as
\[ x_{t+1}=\Pi_{\mathcal{X}}\left( x_t+\frac{\eta}{\sqrt{\sum_{s=1}^t\left\| u_s \right\|_2^2}}u_t \right),  \]
with convention $0/0=0$. The corresponding regret bound is recalled in Appendix~\ref{sec:mirror-descent-type}.
For a given operator $F:\mathcal{X}\to \mathbb{R}^d$, the conversion scheme suggests using the above iteration
with payoff vectors $u_t=F(x_t)-x_t$ (for $t\geqslant 1$), in other words
\[ x_{t+1}=\Pi_{\mathcal{X}}\left( x_t+\frac{\eta}{\sqrt{\sum_{s=1}^t\left\| F(x_s)-x_s \right\|_2^2}}(F(x_t)-x_t) \right).  \]
As stated in Theorem~\ref{thm:adagrad-norm-fp} below, the domain
$\mathcal{X}$ should be known to contain a fixed point $x_*$ and operator
$F$ should be defined at each point of $\mathcal{X}$. 
For $t\geqslant 1$, denote
$\eta_t=\eta/\sqrt{\sum_{s=1}^t\left\| F(x_s)-x_s \right\|_2^2}$. We can see that
above iterations are special cases of KM
iterations \emph{if} $F$ takes values in
$\mathcal{X}$ and \emph{if} $\eta_t<1$ for all $t\geqslant 1$, which is equivalent to
simply having $\eta<\left\| F(x_1)-x_1 \right\|_2$, because then
\[ x_t+\eta_t(F(x_t)-x_t)=(1-\eta_t)x_t+\eta_tF(x_t), \] is a convex
combination of points in $\mathcal{X}$ for all $t\geqslant 1$, and the projection step
$\Pi_\mathcal{X}$ has no effect.
\begin{theorem}[AdaGrad-Norm for fixed points]
  \label{thm:adagrad-norm-fp}
Let \(F:\mathcal{X}\to \mathbb{R}^d\), \(x_*\in \mathcal{X}\) a fixed point of \(F\), \(\eta>0\),
\(x_1\in \mathcal{X}\) that is not a fixed point, and for \(t\geqslant 1\),
\[ x_{t+1}=\Pi_{\mathcal{X}}\left( x_t+\frac{\eta}{\sqrt{\sum_{s=1}^t\left\| F(x_s)-x_s \right\|_2^2}}(F(x_t)-x_t) \right).  \]
Let $T\geqslant 1$ be an integer and assume that
  \[ L_T:=\inf\left\{ L>0,\ \forall 1\leqslant t\leqslant T,\ \left\| F(x_t)-x_t \right\|_2^2\leqslant 2L\left< F(x_t)-x_t, x_*-x_t \right> \right\}  \]
  is finite.
\begin{enumerate}[(i)]
\item\label{item:adagrad-norm-fp-general-bound} Then,
  \begin{multline*}
  \min_{1\leqslant t\leqslant T}\left\| F(x_t)-x_t \right\|_2
  \leqslant \frac{L_T}{\sqrt{T}}\Biggl(
  \frac{\left\| x_1-x_* \right\|_2^2}{\eta}+3\eta+2\eta\log_+ \left( \frac{\eta L_T}{\left\| F(x_1)-x_1 \right\|_2} \right)
  \Biggr),
  \end{multline*}
  where $\log_+=0\vee \log $.
\item\label{item:adagrad-norm-fp-D} Moreover, let 
  \(D\geqslant (\operatorname{diam}_2\mathcal{X})\). Then choosing $\eta=D/\sqrt{2}$ yields
\[ \min_{1\leqslant t\leqslant T}\left\| F(x_t)-x_t \right\|_2\leqslant\frac{2\sqrt{2}DL_T}{\sqrt{T}}. \]
\item\label{item:adagrad-norm-fp-nonexp} Let $L>0$. If $F$ is
  $L$-nonexpansive, then the above holds with \(L_T\leqslant L\) and
  $(x_t)_{t\geqslant 1}$ converges to a fixed point of $F$.
\end{enumerate}
\end{theorem}
\begin{proof}
  First note that $L_T>0$ necessarily because if 
   $L_T=0$, the definition of $L_T$ implies that
  $F(x_t)=x_t$ for all $1\leqslant t\leqslant T$ which contradicts the assumption that
  $x_1$ is not a fixed point.

Using Lemma~\ref{lm:online-to-fp-local} and applying the regret bound for
AdaGrad-Norm from Proposition~\ref{prop:adagrad-norm-regret-bound}
with payoff vectors $u_t=F(x_t)-x_t$ (for $t\geqslant 1$) yields
\begin{align*}
\sum_{t=1}^T\left\| F(x_t)-x_t \right\|_2^2&\leqslant 2L_T\sum_{t=1}^T\left< F(x_t)-x_t, x_{*}-x_t \right> \\
  &\leqslant 2L_T\left( \frac{\max_{1\leqslant t\leqslant T}\left\| x_t-x_* \right\|_2^2}{2\eta}+\eta \right)\sqrt{\sum_{t=1}^T\left\| F(x_t)-x_t \right\|_2^2}. 
\end{align*}
Because $F(x_1)\neq x_1$ by assumption, 
 $\sqrt{\sum_{t=1}^T\left\| F(x_t)-x_t \right\|_2^2}$ is positive and dividing the above by
\[ \sqrt{T\sum_{t=1}^T\left\| F(x_t)-x_t \right\|_2^2} \]
gives
\begin{align}
  \label{eq:2}
  \begin{split}
\min_{1\leqslant t\leqslant T}\left\| F(x_t)-x_t \right\|_2&\leqslant \sqrt{\frac{1}{T}\sum_{t=1}^T\left\| F(x_t)-x_t \right\|_2^2}\\
  &\leqslant \frac{L_T}{\sqrt{T}}\left( \frac{\max_{1\leqslant t\leqslant T}\left\| x_t-x_* \right\|_2^2}{\eta}+2\eta \right).
  \end{split}
\end{align}
In the case where $\max_{1\leqslant t\leqslant T}\left\| x_t-x_* \right\|_2\leqslant D$, the
above bound together with the choice $\eta=D/\sqrt{2}$ yields upper bound~\eqref{item:adagrad-norm-fp-D}.

Let us now prove guarantee~\eqref{item:adagrad-norm-fp-general-bound}. Let $1\leqslant t\leqslant T$. We denote
\[ \eta_t=\frac{\eta}{\sqrt{\sum_{s=1}^t\left\| F(x_s)-x_s \right\|_2^2}} \]
so that
\[ x_{t+1}=\Pi_{\mathcal{X}}\left(x_t+\eta_t(F(x_t)-x_t)  \right). \]
Because $x_*\in \mathcal{X}$, we use the nonexpansiveness of projection
operator $\Pi_{\mathcal{X}}$ to write
\begin{align}
  \label{eq:4}
  \begin{split}
\left\| x_{t+1}-x_* \right\|^2_2&=\left\| \Pi_{\mathcal{X}}(x_t+\eta_t(F(x_t)-x_t)) - \Pi_{\mathcal{X}}(x_*) \right\|_2^2\\
  &\leqslant \left\| x_t+\eta_t(F(x_t)-x_t)-x_* \right\|_2^2\\
  &=\left\| x_t-x_* \right\|_2^2-2\eta_t\left<F(x_t)-x_t , x_*-x_t \right> + \eta_t^2\left\| F(x_t)-x_t \right\|_2^2\\
  &\leqslant \left\| x_t-x_* \right\|_2^2-\frac{\eta_t}{L_T}\left\| F(x_t)-x_t \right\|_2^2+\eta_t^2\left\| F(x_t)-x_t \right\|_2^2\\
  &=\left\| x_t-x_* \right\|_2^2+\left(\eta_t^2-\frac{\eta_t}{L_T} \right)\left\| F(x_t)-x_t \right\|_2^2,
  \end{split}
\end{align}
where the last inequality simply
follows from the definition of $L_T$. Summing, we obtain
\[\left\| x_t-x_* \right\|_2^2\leqslant \left\| x_1-x_* \right\|_2^2+\sum_{s=1}^{t-1}\left(\eta_s^2-\frac{\eta_s}{L_T} \right)\left\| F(x_s)-x_s \right\|_2^2. \]
Let $\tau$ be the largest integer in $\left\{ 1,\dots,T-1 \right\}$ such
that $\eta_{\tau}\geqslant 1/L_T$, or $\tau =0$ if such an integer does not exist, in
which case the above sum only has nonpositive terms and
$\left\| x_t-x_* \right\|_2\leqslant \left\| x_1-x_* \right\|_2$. If $\tau\geqslant 1$, 
because $(\eta_t)_{t\geqslant 1}$ is nonincreasing by definition, for all
$1\leqslant t\leqslant \tau$, it holds that $\eta_t\geqslant 1/L_T$. Then,
\begin{align*}
\left\| x_t-x_* \right\|_2^2&\leqslant \left\| x_1-x_* \right\|_2^2+\sum_{s=1}^{\tau}\left(\eta_s^2-\frac{\eta_s}{L_T} \right)\left\| F(x_s)-x_s \right\|_2^2\\
  &\leqslant \left\| x_1-x_* \right\|_2^2+\sum_{s=1}^{\tau}\eta_s^2\left\| F(x_s)-x_s \right\|_2^2\\
  &=\left\| x_1-x_* \right\|_2^2+\eta^2\left(1+\sum_{s=2}^{\tau}\frac{\eta_s^{-2}-\eta_{s-1}^{-2}}{\eta_s^{-2}}  \right) \\
  &\leqslant \left\| x_1-x_* \right\|_2^2+\eta^2\left( 1+\sum_{s=2}^{\tau}\int_{\eta_{s-1}^{-2}}^{\eta_s^{-2}}\frac{dv}{v} \right) \\
  &=\left\| x_1-x_* \right\|_2^2+\eta^2\left( 1+\log \left( \frac{\eta_{\tau}^{-2}}{\eta_1^{-2}} \right)  \right) \\
  &\leqslant \left\| x_1-x_* \right\|_2^2+\eta^2\left( 1+2\log \left( \frac{\eta L_T}{\left\| F(x_1)-x_1 \right\|_2 } \right) \right) . 
\end{align*}
Plugging the above into~\eqref{eq:2} gives guarantee~\eqref{item:adagrad-norm-fp-general-bound}.

Now consider case~\eqref{item:adagrad-norm-fp-nonexp} where $F$ is
$L$-nonexpansive. $L_T\leqslant L$ follows from
Lemma~\ref{lm:online-to-fp-local}. This makes the above upper bound
finite. Consequently, $(\left\| x_t-x_* \right\|_2^2)_{t\geqslant 1}$ is bounded
and so is $(x_t)_{t\geqslant 1}$. Therefore, there exists $x_{**}$ a
subsequential limit of $(x_t)_{t\geqslant 1}$, which belongs to
$\mathcal{X}$ by closedness of the latter. Combining the second inequality in
\eqref{eq:2} with the above inequality proves that
$\sum_{t=1}^{+\infty}\left\| F(x_t)-x_t \right\|_2^2<+\infty$. Consequently,
$\left\| F(x_t)-x_t \right\|_2^2$ vanishes and, $F$ being continuous
(because $F=(1-L)I+LF_L$ and $F_L$ is nonexpansive), $F(x_{**})=x_{**}$, in other words
$x_{**}$ is a fixed point. Because $(\eta_t)_{t\geqslant 1}$ is positive and
nonincreasing, inequality \eqref{eq:4} implies
\begin{equation}
\label{eq:7}
\left\| x_{t+1}-x_* \right\|_2^2\leqslant \left\| x_t-x_* \right\|_2^2+\eta_1^2\left\| F(x_t)-x_t \right\|_2^2,\quad t\geqslant 1, 
\end{equation}
which is a property called quasi-Fejér monotonicity \citep{ermol1973random}. The following sequence
\[ v_t=\left\| x_t-x_* \right\|_2^2+\eta_1^2\sum_{s=t}^{+\infty}\left\| F(x_s)-x_s \right\|_2^2,\quad t\geqslant 1, \]
is nonnegative and nonincreasing as a
consequence of inequality \eqref{eq:7}. This sequence thus admits a
limit, and thus so does $(\left\| x_t-x_* \right\|_2^2)_{t\geqslant 1}$ because
the other term, being the tail sum of a converging series, vanishes. Note
that this reasoning is true for every fixed point $x_*$ of $F$. It is
true in particular for $x_*=x_{**}$. And the limit of
$\left\| x_t-x_{**} \right\|_2^2$ can only be zero because $x_{**}$ is a
subsequential limit of $(x_t)_{t\geqslant 1}$ by definition. In conclusion,
$(x_t)_{t\geqslant 1}$ converges to $x_{**}$, which is indeed a fixed point.
\end{proof}

The above guarantee is a best-iterate residual bound instead of a
last-iterate one as for KM in Corollary~\ref{cor:km}. Nevertheless, the value of
this minimum and the point achieving it are easy to track, as one only
needs to store the value of the minimum and the corresponding point:
then, at each step, the residual at the new iterate is
computed, compared with the stored minimum, and substituted if lower.

In the case where $F$ is $L$-nonexpansive as in
\eqref{item:adagrad-norm-fp-nonexp} (or more generally, if
$(L_T)_{T\geqslant 1}$ is bounded), the upper bound vanishes at speed
$1/\sqrt{T}$ (similarly to KM iterations) and has the same dependence
on $L$ (up to logarithmic factors) as the bound that would be guaranteed by (projected) KM iterates defined
with operator $F_L:=I+\frac{1}{L}(F-I)$. But running (projected) KM iterations
with operator $F_L$ requires prior knowledge of $L$. The guarantee for
AdaGrad-Norm iterates has the important advantage of being
\emph{adaptive} to unknown $L$.

The upper bound also has a dependence on $\left\| x_1-x_* \right\|_2$,
namely the distance to a solution.
 Without any prior knowledge on this
quantity, the upper bounds grow as the square of this distance. If an
upper bound $D_0$ is known, choosing parameter
$\eta=D_0$ makes the bound depend linearly in $D_0$ (up to logarithmic factors). In the case where
the domain $\mathcal{X}$ is bounded, its diameter can play the role of
$D_0$, as presented in~\eqref{item:adagrad-norm-fp-D} above. An
interesting question is whether it is possible to define iterations
that achieve linear dependence on $\left\| x_1-x_* \right\|_2$ with no
prior knowledge, and is left for future research.

This adaptive guarantee is somewhat analogous to the adaptivity of
AdaGrad-Norm in smooth convex optimization established by \citet{levy2018online}.
In the latter setting however, the guarantee is an upper bound on the
suboptimality gap which is defined in terms of the objective function.
This quantity cannot be expressed nor bounded from above by
the corresponding fixed-point residual alone. Therefore, the above guarantee
for fixed points is not a direct extension of the result for smooth
convex optimization.

\subsection{Adaptivity to positive definite scalings}
\label{sec:adapt-posit-defin}

We now apply the conversion scheme to the AdaGrad-Diagonal and
AdaGrad-Full variants and obtain iterations that are adaptive to
unknown $A$-nonexpansiveness (where $A$ is restricted to positive
diagonal matrices in the former case).

If an operator $F:\mathcal{X}\to \mathcal{X}$, that admits some fixed point
$x_*\in \mathcal{X}$, is $A$-nonexpansive for some known matrix
$A\in \mathcal{S}_{++}^d(\mathbb{R})$, then operator
$F_A:=I+A^{-1}(F-I)$ is nonexpansive with respect to
$\left\|\,\cdot\,\right\|_A$ by definition, and could be used in Projected KM
iterations that would guarantee (replacing
$\left< \,\cdot\,, \,\cdot\, \right>$ by
$\left< \,\cdot\,, \,\cdot\, \right>_A$ in 
Theorem~\ref{thm:proj-KM}):
\[ \min_{1\leqslant t\leqslant T}\left\| F(x_t)-x_t \right\|_{A^{-1}}=\min_{1\leqslant t\leqslant T}\left\| F_A(x_t)-x_t \right\|_A\leqslant \frac{2\left\| x_1-x_* \right\|_A}{\sqrt{T}},\quad T\geqslant 1. \]
The above guarantee depends on $A$ in two places. First, the residual
is measured with $\left\|\,\cdot\,\right\|_{A^{-1}}$, so for a given bound,
the smaller $A$ (in the Loewner order) is, the stronger the guarantee is. Secondly, the upper
bound increases with $A$; so again, the smaller $A$ is, the stronger
the guarantee. In the case where an estimate on the diameter of the
domain $\mathcal{X}$ is available, iterates based on AdaGrad-Diagonal and
AdaGrad-Full obtain in Theorems~\ref{thm:adagrad-diagonal-fp} and
\ref{thm:adagrad-full-fp} below, without prior knowledge of $A$, somewhat
similar guarantees which are upper bounds on the best-iterate residual
$\min_{1\leqslant t\leqslant T}\left\| F_A(x_t)-x_t \right\|_A$ that have the usual
$1/\sqrt{T}$ dependence and that grow with $A$, specifically as
$\sqrt{\operatorname{Tr}A}$.

For a vector $u\in \mathbb{R}^d$, $\operatorname{diag}u$ denotes the diagonal
matrix of size $d$ with diagonal coefficients given by $u$. For a
matrix $A\in \mathcal{S}_{++}^d(\mathbb{R})$, $A^{-1/2}$ denotes the inverse of the
square-root of $A$ and the projection onto $\mathcal{X}$ with respect to
Mahalanobis norm $\left\|\,\cdot\,\right\|_A$ is denoted $\Pi_{\mathcal{X},A}$.

\begin{theorem}[Guarantee for AdaGrad-Diagonal for fixed points]
  \label{thm:adagrad-diagonal-fp}
Let \(F:\mathcal{X}\to \mathbb{R}^d\), \(x_*\in \mathcal{X}\) a fixed point of \(F\), \(\eta,\varepsilon>0\), \(x_1\in \mathcal{X}\), and for \(t\geqslant 1\),
\[ x_{t+1}=\Pi_{\mathcal{X},A_t}(x_t+A_t^{-1}(F(x_t)-x_t)),  \]
where
\[ A_t=\eta^{-1}\operatorname{diag}\left( \sqrt{\varepsilon^2+\sum_{s=1}^t(F(x_s)_i-x_{s,i})^2} \right)_{1\leqslant i\leqslant d}.   \]
Let \(T\geqslant 1\) be an integer and let
\(A\in\mathcal{S}_{++}^d(\mathbb{R})\) be a diagonal matrix.
\begin{enumerate}[(i)]
\item\label{item:adagrad-diagonal-regret-fp} If, for all $1\leqslant t\leqslant T$,
  \[ \left\| F(x_t)-x_t \right\|_{A^{-1}}^2\leqslant 2\left< F(x_t)-x_t, x_*-x_t \right>, \]
  then
  \[ \min_{1\leqslant t\leqslant T}\left\| F(x_t)-x_t \right\|_{A^{-1}}\leqslant \sqrt{\frac{\operatorname{Tr}A}{T}}\left( \frac{D_{\infty,T}^2}{\eta}+2\eta \right)+\frac{\varepsilon d}{\sqrt{T(\operatorname{Tr}A)}}, \]
  where $D_{\infty,T}=\max_{1\leqslant t\leqslant T}\left\| x_t-x_* \right\|_{\infty}$ and consequently
  \[ \min_{1\leqslant t\leqslant T}\left\| F(x_t)-x_t \right\|_2\leqslant \sqrt{\frac{(\operatorname{Tr}A)\lambda_{\mathrm{max}}(A)}{T}}\left( \frac{D_{\infty,T}^2}{\eta}+2\eta \right)+\frac{\varepsilon d}{\sqrt{T}} . \]
\item\label{item:adagrad-diagonal-radius-fp} Moreover, let
  \(D\geqslant \operatorname{diam}_{\infty}(\mathcal{X})\). Then choosing
  \(\eta=D/\sqrt{2}\) yields, 
  \[ \min_{1\leqslant t\leqslant T}\left\| F(x_t)-x_t \right\|_{A^{-1}}\leqslant 2\sqrt{2}D\sqrt{\frac{\operatorname{Tr}A}{T}}+\frac{\varepsilon d}{\sqrt{T(\operatorname{Tr}A)}} \]
  and
  \[ \min_{1\leqslant t\leqslant T}\left\| F(x_t)-x_t \right\|_{2}\leqslant 2\sqrt{2}D\sqrt{\frac{(\operatorname{Tr}A)\lambda_{\mathrm{max}}(A)}{T}}+\frac{\varepsilon d}{\sqrt{T}}.  \]
  \item\label{item:adagrad-diagonal-moreover} In particular, if operator $F$ is $A$-nonexpansive, the above guarantees~\eqref{item:adagrad-diagonal-regret-fp} and~\eqref{item:adagrad-diagonal-radius-fp} hold.
\item\label{item:adagrad-diagonal-cv} If operator $F$ is
  $A$-nonexpansive and $(x_t)_{t\geqslant 1}$ is bounded, then it converges to
  a fixed point of $F$.
\end{enumerate}
\end{theorem}
\begin{proof}
  If $x_1=x_*$, the iteration gives $x_t=x_*$ for all $t\geqslant 1$ and the
  results are true. From now on, we assume $x_1\neq x_*$.
  Summing inequality
  \begin{equation}
\label{eq:9}
\left\|F(x_t)-x_t\right\|_{A^{-1}}^2\leqslant 2\left< F(x_t)-x_t, x_*-x_t \right>,\quad 1\leqslant t\leqslant T,  
\end{equation}
  and applying the regret bound for
  AdaGrad-Diagonal from Proposition~\ref{prop:adagrad-diagonal} gives
  \begin{multline*}
    \sum_{t=1}^T\left\| F(x_t)-x_t \right\|_{A^{-1}}^2\\\leqslant 2\left( \frac{\varepsilon}{2\eta}\left\| x_1-x_* \right\|^2_2+\left( \frac{D_{\infty,T}^2}{2\eta}+\eta \right)\sqrt{(\operatorname{Tr}A)\sum_{t=1}^T\left\| F(x_t)-x_t \right\|_{A^{-1}}^2}  \right), 
  \end{multline*}
  where $D_{\infty,T}=\max_{1\leqslant t\leqslant T}\left\| x_t-x_* \right\|_{\infty}$. The above
  shows that quantity
  $X=\sqrt{\sum_{t=1}^T\left\| F(x_t)-x_t \right\|_{A^{-1}}^2}$ satisfies
  inequality $X^2\leqslant \alpha+\beta X$, where
  \[ \alpha=\frac{\varepsilon}{\eta}\left\| x_1-x_* \right\|_2^2\quad \operatorname{and}\quad \beta =\left( \frac{D_{\infty,T}^2}{\eta}+2\eta \right)\sqrt{\operatorname{Tr}A}. \]
  We deduce that 
  \begin{align*}
  X&\leqslant\frac{\beta+\sqrt{\beta^2+4\alpha}}{2}\leqslant \frac{1}{2}\left( \beta+\beta\sqrt{1+\frac{4\alpha}{\beta^2}} \right)\leqslant \frac{1}{2}\left( \beta+\beta\left( 1+\frac{2\alpha}{\beta^2} \right)  \right)   \\
    &=\beta+\frac{\alpha}{\beta}=\beta+\frac{\varepsilon(\operatorname{Tr}A)^{-1/2}\left\| x_1-x_* \right\|_2^2}{\max_{1\leqslant t\leqslant T}\left\| x_t-x_* \right\|_{\infty}^2+2\eta^2}\\
    &\leqslant \beta+\frac{\varepsilon(\operatorname{Tr}A)^{-1/2}\left\| x_1-x_* \right\|_2^2}{\left\| x_1-x_* \right\|_{\infty}^2}\leqslant \beta+\frac{\varepsilon d}{\sqrt{\operatorname{Tr}A}},
  \end{align*}
  where we could divide by $\left\| x_1-x_* \right\|_{\infty}^2$ because
  $x_1\neq x_*$ as assumed above.
  Then, dividing by $\sqrt{T}$, we obtain
  \begin{align*}
\min_{1\leqslant t\leqslant T}\left\| F(x_t)-x_t \right\|_{A^{-1}}&\leqslant \sqrt{\frac{1}{T}\sum_{t=1}^T\left\| F(x_t)-x_t \right\|_{A^{-1}}^2}\\
    &\leqslant \sqrt{\frac{\operatorname{Tr}A}{T}}\left( \frac{D_{\infty,T}^2}{\eta}+2\eta \right)+ \frac{\varepsilon d}{\sqrt{T (\operatorname{Tr}A)}},
\end{align*}
which gives the first inequality from
\eqref{item:adagrad-diagonal-regret-fp}. The other inequality on
$\min_{1\leqslant t\leqslant T}\left\| F(x_t)-x_t \right\|_2$ is simply obtained by
using inequalities
\[ \left\|\,\cdot\,\right\|_2\leqslant \sqrt{\lambda_{\mathrm{max}}(A)}\left\|\,\cdot\,\right\|_{A^{-1}}\quad \text{and}\quad \lambda_{\mathrm{max}}(A)\leqslant \operatorname{Tr}A. \]
\eqref{item:adagrad-diagonal-radius-fp} follows
immediately.
Case~\eqref{item:adagrad-diagonal-moreover} simply corresponds to the
sufficient condition given by Lemma~\ref{lm:online-to-fp}.

\eqref{item:adagrad-diagonal-cv}
We prove that the iterates satisfy a variable-metric quasi-Fejér property~\citep{combettes2013variable}.
For all $t\geqslant 1$, write $u_t=F(x_t)-x_t$.
$F$ being $A$-nonexpansive, inequality \eqref{eq:9} is satisfied
for all $t\geqslant 1$. The regret bound above gives
\[
 S_T:=\sum_{t=1}^T\left\|u_t\right\|_{A^{-1}}^2 \leqslant c_0+c_1\sqrt{S_T}
\]
for all $T\geqslant 1$, where $c_0,c_1$ are finite nonnegative constants independent of $T$.
We deduce that $\sum_{t=1}^{+\infty}\left\|u_t\right\|_2^2<+\infty$.
 By this summability, each diagonal entry
\[
 (A_t)_{ii}=\frac1\eta\sqrt{\varepsilon^2+\sum_{s=1}^t u_{s,i}^2}
\]
is nondecreasing and converges to the finite value
$\eta^{-1}\sqrt{\varepsilon^2+\sum_{s=1}^{+\infty}u_{s,i}^2}$.
Thus $A_t\to A_\infty$, with
\[
 \frac{\varepsilon}{\eta}I_d\preceq A_t\preceq A_{t+1} \]
and thus
\[ \sum_{t=1}^{+\infty}\operatorname{Tr}(A_{t+1}-A_t) =\operatorname{Tr}(A_\infty-A_1)<+\infty. \]

Let 
$R:=\sup_{t\geqslant 1}\left\|x_t-x_*\right\|_2$ which is finite because
$(x_t)_{t\geqslant 1}$ is here bounded by assumption.
For $t\geqslant 1$, from the proof of Proposition~\ref{prop:omd} we obtain
\[ \left\| x_{t+1}-x_* \right\|_{A_t}^2\leqslant\left\| x_t-x_* \right\|_{A_t}^2+2\left< u_t, x_t-x_* \right> + \left\| u_t \right\|_{A_t^{-1}}^2,  \]
which, combined with \eqref{eq:9} yields
\[ \left\|x_{t+1}-x_*\right\|_{A_{t+1}}^2 \leqslant \left\|x_t-x_*\right\|_{A_t}^2+e_t, \]
where
$e_t=(\eta/\varepsilon)\left\|u_t\right\|_2^2
+R^2\operatorname{Tr}(A_{t+1}-A_t)$ is summable.
Consequently,
$(\left\|x_t-x_*\right\|_{A_t}^2+\sum_{s=t}^{+\infty}e_s)_{t\geqslant 1}$
is nonnegative and nonincreasing, so
$(\left\|x_t-x_*\right\|_{A_t}^2)_{t\geqslant 1}$ has a limit.
The argument applies to every fixed point $x_*$.
Boundedness of the iterates gives the
existence of a subsequential limit $x_{**}$, which is a fixed point by continuity of $F$ and $u_t\to 0$.
Taking $x_*=x_{**}$, boundedness of $(A_t)_{t\geqslant 1}$ implies that the
limit of $(\left\|x_t-x_{**}\right\|_{A_t}^2)_{t\geqslant 1}$ can only be zero because
$x_{**}$ is a subsequential limit of $(x_t)_{t\geqslant 1}$ by definition.
Together with $A_t\succeq(\varepsilon/\eta)I_d$, this implies $x_t\to x_{**}$.
\end{proof}

\begin{theorem}[Guarantee for AdaGrad-Full for fixed points]
  \label{thm:adagrad-full-fp}
Let \(F:\mathcal{X}\to \mathbb{R}^d\), \(x_*\in \mathcal{X}\) a fixed point of \(F\), \(\eta,\varepsilon>0\), \(x_1\in \mathcal{X}\), and for \(t\geqslant 1\),
\[ x_{t+1}=\Pi_{\mathcal{X},A_t}(x_t+A_t^{-1}(F(x_t)-x_t)),  \]
where
\[ A_t=\eta^{-1}\left( \varepsilon^2I_d+\sum_{s=1}^t(F(x_s)-x_s)(F(x_s)-x_s)^{\top} \right)^{1/2}.   \]
Let \(T\geqslant 1\) be an integer and $A\in \mathcal{S}_{++}^d(\mathbb{R})$. 
\begin{enumerate}[(i)]
\item\label{item:adagrad-full-fp} If for all $1\leqslant t\leqslant T$,
  \[ \left\| F(x_t)-x_t \right\|_{A^{-1}}^2\leqslant 2\left< F(x_t)-x_t, x_*-x_t \right> , \]
  then
  \[ \min_{1\leqslant t\leqslant T}\left\| F(x_t)-x_t \right\|_{A^{-1}}\leqslant \sqrt{\frac{\operatorname{Tr}A}{T}}\left( \frac{D_{2,T}^2}{\eta}+2\eta \right)+\frac{\varepsilon}{\sqrt{T(\operatorname{Tr}A)}}, \]
  where $D_{2,T}=\max_{1\leqslant t\leqslant T}\left\| x_t-x_* \right\|_2$ and consequently,
  \[ \min_{1\leqslant t\leqslant T}\left\| F(x_t)-x_t \right\|_{2}\leqslant \sqrt{\frac{\lambda_{\mathrm{max}}(A)(\operatorname{Tr}A)}{T}}\left( \frac{D_{2,T}^2}{\eta}+2\eta \right)+\frac{\varepsilon}{\sqrt{T}}, \]
\item\label{item:adagrad-full-fp-D} Moreover, let
  \(D\geqslant \operatorname{diam}_2(\mathcal{X})\). Then, choosing \(\eta=D/\sqrt{2}\) yields,
  \[ \min_{1\leqslant t\leqslant T}\left\| F(x_t)-x_t \right\|_{A^{-1}}\leqslant 2D\sqrt{\frac{2\operatorname{Tr}A}{T}}+\frac{\varepsilon}{\sqrt{T(\operatorname{Tr}A)}}, \]
  and therefore
  \[ \min_{1\leqslant t\leqslant T}\left\| F(x_t)-x_t \right\|_{2}\leqslant 2D\sqrt{\frac{2\lambda_{\mathrm{max}}(A)(\operatorname{Tr}A)}{T}}+\frac{\varepsilon}{\sqrt{T}}. \]
\item\label{item:adagrad-full-fp-nonexp} In particular, if $F$ is
  $A$-nonexpansive, the above guarantees \eqref{item:adagrad-full-fp}
  and \eqref{item:adagrad-full-fp-D} hold.
\item\label{item:adagrad-full-cv} If operator $F$ is
  $A$-nonexpansive and $(x_t)_{t\geqslant 1}$ is bounded, then it converges to
  a fixed point of $F$.
\end{enumerate}
\end{theorem}
\begin{proof}
  The residual bounds follow as in Theorem~\ref{thm:adagrad-diagonal-fp},
  using Proposition~\ref{prop:adagrad-full}. The convergence
  argument for~\eqref{item:adagrad-full-cv} is similar to that in
  Theorem~\ref{thm:adagrad-diagonal-fp}. Here, writing $u_t=F(x_t)-x_t$,
  the series $\sum_s u_su_s^\top$ converges absolutely in operator norm,
  since $\|u_su_s^\top\|_{\mathrm{op}}=\|u_s\|_2^2$ is summable.
  Continuity and monotonicity of the matrix square root therefore give
  $A_t\to A_\infty$ and
  $(\varepsilon/\eta)I_d\preceq A_t\preceq A_{t+1}$, as required.
\end{proof}

In both iterations, the only role of $\varepsilon>0$ is to make $A_t$
invertible, and is taken very small in practice, e.g.,\ $10^{-10}$.
In the unconstrained case $\mathcal{X}=\mathbb{R}^d$, both AdaGrad-Diagonal and
AdaGrad-Full can be simply written
\[ x_{t+1}=x_t+A_t^{-1}(F(x_t)-x_t),\quad t\geqslant 1. \]

The assumption relating operator $F$ and matrix $A$ in
\eqref{item:adagrad-full-fp} in Theorems~\ref{thm:adagrad-diagonal-fp}
and \ref{thm:adagrad-full-fp} is, as stated in
\eqref{item:adagrad-full-fp-nonexp}, a relaxation of $A$-nonexpansiveness
 focusing on solution $x_*$, and which is local along the trajectory of iterates $x_1,\dots,x_T$.

The main guarantees in Theorems~\ref{thm:adagrad-diagonal-fp}
and \ref{thm:adagrad-full-fp} are of \emph{best-iterate} type: more
precisely, they give upper bounds on
\[ \min_{1\leqslant t\leqslant T}\left\| F(x_t)-x_t \right\|_{A^{-1}}, \]
meaning the lowest
fixed-point residual so far, measured with norm
$\left\|\,\cdot\,\right\|_{A^{-1}}$, where $A$ satisfies the assumption from
\eqref{item:adagrad-full-fp}. This quantity is not tractable as
matrix $A$ may be unknown. In particular, one does not know which
point among $x_1,\dots,x_T$ actually enjoys the upper bound on its
fixed-point residual. Fortunately, a guarantee on 
$\min_{1\leqslant t\leqslant T}\left\| F(x_t)-x_t \right\|_{2}$ is immediately deduced
at the cost of a factor $\sqrt{\lambda_{\text{max}}(A)}$ in the bound.
Although this is a weaker guarantee, the quantity is at least tractable.

Regarding AdaGrad-Diagonal, the adaptivity to the best-in-hindsight 
positive diagonal preconditioner achieved in
Theorem~\ref{thm:adagrad-diagonal-fp} is analogous to the results
obtained for AdaGrad-Diagonal in the context of smooth convex
optimization \citep{liu2024large,jiang2024convergence}, where the
corresponding property is called \emph{anisotropic smoothness} or
\emph{coordinate-wise smoothness}.

In the case of an unbounded domain $\mathcal{X}$, unlike the guarantee for
AdaGrad-Norm from Theorem~\ref{thm:adagrad-norm-fp}, the above
statements do not always guarantee a $1/\sqrt{T}$ speed of convergence
over fixed-point residuals, because the distance to solutions
($D_{\infty,T}$ and $D_{2,T}$ in Theorems~\ref{thm:adagrad-diagonal-fp} and
\ref{thm:adagrad-full-fp}, respectively) appear in the upper bound and
are not known to be bounded in general. In the context of
unconstrained smooth convex optimization however, AdaGrad-Diagonal has
been proved to ensure the convergence of the iterates
\citep{traore2021sequential}, and therefore their boundedness. Whether
such a guarantee can be adapted to fixed points in the unconstrained
or unbounded case remains an open question. In the case where the
domain \emph{is} bounded, or if the iterates are known to be bounded
for some reason, we recover a $1/\sqrt{T}$ convergence rate. We now
restrict the discussion to that case and focus on the dependence of
the upper bounds on the diameter of the domain and on the matrix $A$.

For this discussion, we neglect the term
in $\varepsilon$, as it is chosen very small in practice.
Regarding the dependence on $A$ and in the diameter of the domain $\mathcal{X}$,
the upper bound for AdaGrad-Diagonal scales as
$\sqrt{\operatorname{Tr}A}(\operatorname{diam}_{\infty}\mathcal{X})$, whereas, as
discussed in the introduction, the KM
iterations run with operator $F_A:=I+A^{-1}(F-I)$ would guarantee a
dependence on $\operatorname{diam}_A\mathcal{X}$, which denotes the diameter
measured with $\left\|\,\cdot\,\right\|_A$. Although the latter quantity
is lower in general, the upper bound for AdaGrad-Diagonal is in $\sqrt{\operatorname{Tr}A}(\operatorname{diam}_{\infty}\mathcal{X})$, without prior knowledge of
$A$, which is remarkable.  The upper bound achieved by
AdaGrad-Full is seemingly weaker, as it scales as
$\sqrt{\operatorname{Tr}A}(\operatorname{diam}_2\mathcal{X})$, but it is
adaptive to all positive definite matrices $A$ instead of just
diagonal ones.

Furthermore, let us compare the above guarantees for AdaGrad-Diagonal
and AdaGrad-Full with the guarantee achieved by AdaGrad-Norm in
Theorem~\ref{thm:adagrad-norm-fp}. Direct comparison is difficult,
because for AdaGrad-Diagonal and AdaGrad-Full, 
the assumption is finer and the main guarantee is a bound on the fixed-point residuals measured with $\left\|\,\cdot\,\right\|_{A^{-1}}$ and not
$\left\|\,\cdot\,\right\|_2$ as for AdaGrad-Norm. Of course, bounds on
Euclidean fixed-point residuals can be deduced but are weaker than the
main guarantee. We can nevertheless try the following comparison.
Let us first assume that diagonal matrix $A\in \mathcal{S}_{++}^d(\mathbb{R})$ is such that
$F$ is $A$-nonexpansive. Then, the $L_T$ coefficient
from Theorem~\ref{thm:adagrad-norm-fp} can be bounded as follows:
\begin{align*}
L_T&:=\sup_{1\leqslant t\leqslant T}\frac{\left\| F(x_t)-x_t \right\|_2^2}{2\left< F(x_t)-x_t, x_*-x_t \right> }\\
  &\leqslant \sup_{1\leqslant t\leqslant T}\frac{\lambda_{\text{max}}(A)\left\| F(x_t)-x_t \right\|_{A^{-1}}^2}{2\left< F(x_t)-x_t, x_*-x_t \right> }\\
  &\leqslant \lambda_{\text{max}}(A), 
\end{align*}
where convention $0/0=0$ is used.
Then, AdaGrad-Norm guarantees on the best Euclidean residual
$\min_{1\leqslant t\leqslant T}\left\| F(x_t)-x_t \right\|_2$ an upper bound that
scales as $\operatorname{diam}_2(\mathcal{X})\lambda_{\text{max}}(A)$, whereas
AdaGrad-Diagonal guarantees an upper bound on the best Euclidean
residual (of its own iterates) that
scales as $\operatorname{diam}_{\infty}(\mathcal{X})\sqrt{\lambda_{\text{max}}(A)\operatorname{Tr}A}$.
Comparing these two dependencies, the deduced bound for
AdaGrad-Diagonal is stronger in cases where the ratio
$(\operatorname{diam}_{\infty}\mathcal{X})/(\operatorname{diam}_2\mathcal{X})$ (which can be as
low as $1/\sqrt{d}$, depending on the shape of the domain) is significantly lower
than $\sqrt{\lambda_{\text{max}}(A)/(\operatorname{Tr}A)}$, which can be as
high as $1$. But note that this comparison is biased towards AdaGrad-Norm, as the main
guarantee for AdaGrad-Diagonal is on a different quantity.

Direct comparison is even more difficult with AdaGrad-Full because,
compared to AdaGrad-Diagonal and AdaGrad-Norm it
does not offer a better bound on Euclidean residuals
$\min_{1\leqslant t\leqslant T}\left\| F(x_t)-x_t \right\|_2$, but its advantage only comes
from the fact that its main guarantee is a bound on
$\min_{1\leqslant t\leqslant T}\left\| F(x_t)-x_t \right\|_{A^{-1}}$, which may be a
much higher quantity, because $A$ may be any positive definite matrix
satisfying assumption from \eqref{item:adagrad-full-fp}.

Unlike AdaGrad-Diagonal, we believe that the interest of AdaGrad-Full
is mostly theoretical, because computing $A_t^{-1}(F(x_t)-x_t)$
requires an eigendecomposition and scales as $O(d^3)$. It may however
be interesting to consider a variant that restricts to a class of
block-diagonal matrices that remains scalable.

\section{Numerical experiments}
\label{sec:numer-exper}

We examine the practical advantage of using AdaGrad-based iterations
(as well as a few variants) for solving fixed-point problems. We
consider three problems that are commonly solved by methods that can
be interpreted as fixed-point iterations: finding the stationary
distribution of a Markov chain, image denoising using total-variation
regularization and solving two-player zero-sum games. These
experiments are only presented as illustrations, and contain
comparisons with a given baseline method only. They do not aim at
competing with state-of-the-art methods for each problem, which would
be out of the scope of the present work.

The comparisons concern individual instances and iteration
counts, not average performance over random instances or running time.
Statements about convergence below describe the observed finite
trajectories, rather than establishing asymptotic convergence.

\subsection{AdaGrad variants as heuristics}
\label{sec:heur-adagr-vari}
In addition to AdaGrad-Norm and AdaGrad-Diagonal algorithms, we
consider adaptations of RMSprop~\citep{tieleman2012lecture} and
Adam~\cite{kingma2015adam} algorithms, which are AdaGrad variants that
are extremely efficient in the context of deep learning. Their
theoretical understanding is not as solid as for AdaGrad, and we only
consider them as heuristics. RMSprop replaces the plain sum in the
square-root of the scaling by an exponentially weighted sum. Adam also
replaces the step direction vector $F(x_t)-x_t$ by an exponential
moving sum. We consider unnormalized sums, without the factors
$1-\beta$ and $1-\alpha$ or Adam bias correction.
Let $\varepsilon_0>0$ and $\varepsilon_{0,i}>0$ for each $i\in \{1,\dots,d\}$. The updates are, for
$t\geqslant 1$,

\begin{equation}
\label{eq:rmsprop-norm}
\tag{RMSprop-Norm}
x_{t+1}=\Pi_{\mathcal{X}}\left( x_t+\frac{\eta}{\sqrt{\beta^t \varepsilon_0^2+\sum_{s=1}^t\beta^{t-s}\left\| F(x_s)-x_s \right\|_2^2}}(F(x_t)-x_t) \right).
\end{equation}
\begin{equation}
\label{eq:adam-norm}
\tag{Adam-Norm}
\begin{split}
x_{t+1}&=\Pi_{\mathcal{X}}\Biggl(
x_t+\frac{\eta}{\sqrt{\beta^t \varepsilon_0^2+\sum_{s=1}^t\beta^{t-s}\left\| F(x_s)-x_s \right\|_2^2}}\\
&\qquad {}\times\sum_{s=1}^t\alpha^{t-s}(F(x_s)-x_s)\Biggr).
\end{split}
\end{equation}
\begin{align}
\label{eq:rmsprop-diagonal}
\tag{RMSprop-Diagonal}
 \begin{split}
  A_t&=\eta^{-1}\operatorname{diag}\left( \sqrt{\beta^t \varepsilon_{0,i}^2+\sum_{s=1}^t\beta^{t-s}(F(x_s)_i-x_{s,i})^2} \right)_{1\leqslant i\leqslant d},\\
x_{t+1}&=\Pi_{\mathcal{X},A_t}\left( x_t+A_t^{-1}(F(x_t)-x_t) \right),\quad t\geqslant 1.
 \end{split} 
\end{align}

\begin{align}
\label{eq:adam-diagonal}
\tag{Adam-Diagonal}
 \begin{split}
  A_t&=\eta^{-1}\operatorname{diag}\left( \sqrt{\beta^t \varepsilon_{0,i}^2+\sum_{s=1}^t\beta^{t-s}(F(x_s)_i-x_{s,i})^2} \right)_{1\leqslant i\leqslant d},\\
x_{t+1}&=\Pi_{\mathcal{X},A_t}\left( x_t+A_t^{-1}\left(\sum_{s=1}^t\alpha^{t-s}(F(x_s)-x_s)  \right)  \right),\quad t\geqslant 1.
 \end{split} 
\end{align}
The RMSprop implementations use $\beta=.99$ for the Markov chain
and games, and $\beta=.999$ for image denoising. Adam uses
$\beta=.999$ and $\alpha=.9$ throughout, with zero initial momentum.
The squared accumulators start at zero for the Markov Norm variants
and all game variants, at $10^{-8}$ in each coordinate for the Markov
Diagonal variants, and at $10^{-10}$ for all image-denoising variants.
Thus the initial accumulator decays in RMSprop and Adam, whereas it
remains an additive constant in AdaGrad. No further constant is added
outside the square root. Values for $\eta$ are obtained from tuning.
We do not implement AdaGrad-Full
iterations, as their precise computation does not scale well for large
problems.

\subsection{Markov chain stationary distribution computation}
\label{sec:mark-chain-stat}
We construct a Markov chain that is difficult to explore by splitting
the states into two almost-disconnected clusters. Let $n\geqslant 1$ and
$\left\{ 1,\dots,2n\right\}$ be the set of states. Let
$\left\{ 1,\dots,n \right\}$ and $\left\{ n+1,\dots,2n \right\}$ be
the two clusters. States $1,n,n+1$ and $2n$ are called \emph{boundary states} and
the others are called \emph{interior states}. Let $p\in (0,1)$ be the
probability of jumping from one cluster to the other, uniformly.
Conditionally on not jumping to the other cluster, an interior state
(resp.\ boundary state) has probability $1/3$ (resp.\ $1/2$) to stay
put and $1/3$ to move to each neighbor (resp.\ $1/2$ to move to its
unique neighbor). These probabilities are encoded in the transition
matrix $P\in \mathbb{R}^{2n\times 2n}$.

We aim at iteratively computing the stationary distribution $\pi_*$ that
belongs to $\mathcal{X}=\Delta_{2n}$ (the unit simplex in
$\mathbb{R}^{2n}$). In this section, probability distributions are represented
as row vectors, so the stationary distribution is defined by the fixed-point property
$\pi_* P=\pi_*$. The basic \emph{power iteration} corresponds to
$\pi_{t+1}=\pi_tP$ (for $t\geqslant 1$) \citep{kemeny1960finite}.
We consider $n=5000$ in the notation above and 
$p=10^{-8}$.

\begin{figure}
  \centering
\resizebox{225pt}{!}{
  \input{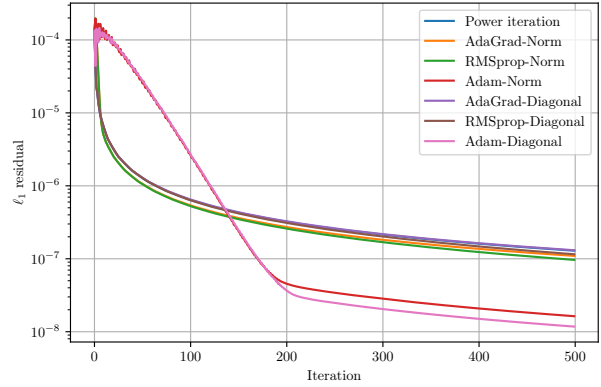}}
\caption{Computation of the stationary distribution of a Markov chain.}\label{fig:markov-chain}
\end{figure}

All methods start from the uniform distribution on the $2n$ states.
The transition matrix is stored densely and its rows are normalized
after assembly. The comparison uses at most $500$ iterations, with
stopping tolerance $10^{-10}$ on $\|\pi_tP-\pi_t\|_1$, measured
before the update. The selected step-sizes in the source are
$\eta=10^{-4}$ for all adaptive methods except RMSprop-Norm, for
which $\eta=3\cdot10^{-5}$. The Norm variants use Euclidean simplex
projection, implemented by sorting, and the Diagonal variants use
weighted simplex projection by bisection, with a mass tolerance of
$10^{-10}$ and at most $100$ bisection steps.
Convergence is observed through the $\ell_1$ residual and is plotted in
Figure~\ref{fig:markov-chain}.

AdaGrad-Norm, RMSprop-Norm, AdaGrad-Diagonal and RMSprop-Diagonal
achieve similar convergence as the power iteration. Adam-Norm and
Adam-Diagonal however, after an initial slower convergence, quickly
reach lower values than the other methods, demonstrating a potential
benefit for this problem.

\subsection{Total variation image denoising with the Chambolle--Pock operator}
\label{sec:cons-optim-with}
Let $m,n\geqslant 1$. We consider the classical ROF image denoising model~\citep{rudin1992nonlinear}:
\[ \min_{u\in \mathbb{R}^{m\times n}}\left\{ \frac{1}{2}\left\| u-f \right\|_2^2+\lambda\left\| \nabla u \right\|_1 \right\},  \]
where $f\in \mathbb{R}^{m\times n}$ denotes a noisy grayscale image of size $m\times n$, $u\in \mathbb{R}^{m\times n}$ the
output (denoised) image, $\nabla u$ the spatial gradient, and $\lambda>0$ a regularization
parameter. The celebrated Chambolle--Pock (CP) algorithm~\citep{chambolle2011first} solves the associated
primal-dual problem and can be written as KM iterations associated
with the following operator, defined for
$(u,p)\in \mathbb{R}^{m\times n}\times \mathbb{R}^{2(m\times n)}$ as
\[ F_{\tau,\sigma,\lambda}(u,p)=\left(2u'-u,2\cdot \Pi_{\left\|\,\cdot\,\right\|_{\infty}\leqslant \lambda}(p+\sigma\nabla \bar{u})-p\right), \]
where
\[ u'=\frac{u+\tau \operatorname{div}p+\tau f}{1+\tau}\quad \text{and}\quad \bar{u}=u'+\theta(u'-u), \]
and where $\operatorname{div}$ denotes the backward divergence
operator and
$\Pi_{\left\|\,\cdot\,\right\|_{\infty}\leqslant \lambda}$ the component-wise projection onto the
$\ell^{\infty}$ ball of radius $\lambda$. AdaGrad-based iterations on
$\mathcal{X}=\mathbb{R}^{m\times n}\times \mathbb{R}^{2(m\times n)}$ and variants are implemented with the above
operator.

The code uses forward differences in both spatial directions, with
zero horizontal differences in the last column and zero vertical
differences in the last row. Thus $\|\nabla u\|_1$ is the
anisotropic total variation. The initialization is $u=f$ and $p=0$;
the corresponding inactive boundary components of $p$ remain zero.
On this subspace, the implemented divergence is the negative adjoint
of the gradient. The CP baseline corresponds to KM iterations with
operator $F_{\tau,\sigma,\lambda}$.

The image we use is the standard boat test
image\footnote{\url{https://sipi.usc.edu/database/preview/misc/boat.512.png}}
of size $512\times 512$. Pixel values are scaled to $[0,1]$, and
additive Gaussian noise with standard deviation $0.1$ is applied
without clipping. We use extrapolation parameter $\theta=1$ and
regularization parameter $\lambda=0.1$.

\begin{figure}[p]
  \centering
  \begin{subfigure}[c]{1.0\linewidth}
  \centering
\resizebox{225pt}{!}{\input{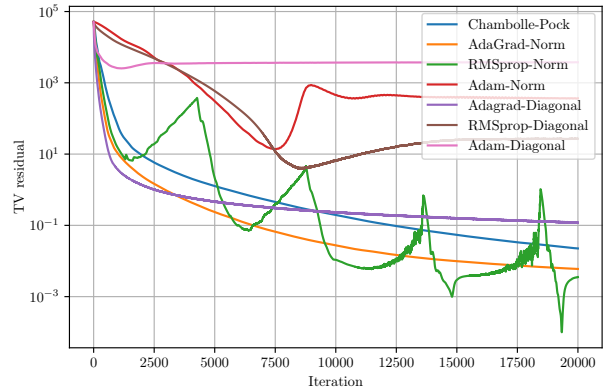}}
\caption{$\tau=\sigma=10^{-2}$, the Chambolle--Pock algorithm converges slowly.}
  \end{subfigure}
  \par\smallskip
  \begin{subfigure}[c]{1.0\linewidth}
  \centering
\resizebox{225pt}{!}{\input{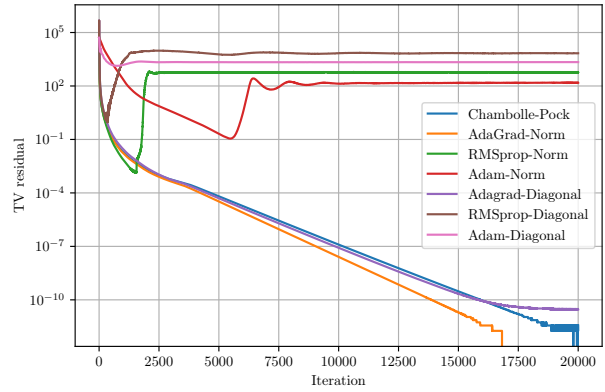}}
\caption{$\tau=\sigma=0.2$, the Chambolle--Pock algorithm converges quickly.}
  \end{subfigure}
  \par\smallskip
  \begin{subfigure}[c]{1.0\linewidth}
  \centering
\resizebox{225pt}{!}{\input{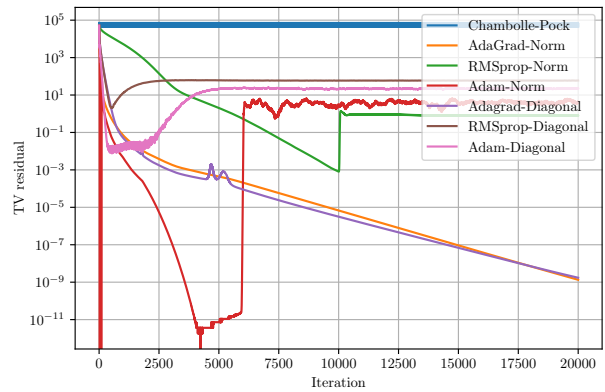}}
\caption{$\tau=\sigma=1$, the Chambolle--Pock algorithm does not converge.}
  \end{subfigure}
\caption{Total-variation image denoising.}\label{fig:image}
\end{figure}

The convergence is plotted in Figure~\ref{fig:image} where the
monitored quantity (labeled ``TV residual'') is $\|\nabla u'\|_1-c$, where $u'$ is computed
inside the operator, before the outer adaptive update, and $c$ is a
minimum value computed with tuned methods and a large number of
iterations. We try three different values for the step-sizes $\tau$ and
$\sigma$, to which the CP algorithm is very sensitive.

For $\tau=\sigma=10^{-2}$, the convergence of the CP algorithm is relatively
slow. The AdaGrad-Norm and AdaGrad-Diagonal algorithms also converge
slowly. The AdaGrad-Norm algorithm seems to converge slightly faster
than CP and the AdaGrad-Diagonal converges faster than CP during the
initial iterations only. The RMSprop-Norm algorithm reaches lower values
faster than the other algorithms, but is unstable, and the overall
convergence speed seems comparable to the other converging algorithms.
Remaining algorithms Adam-Norm, RMSprop-Diagonal and Adam-Diagonal do
not converge.

For $\tau=\sigma=0.2$, the CP algorithm, AdaGrad-Norm and AdaGrad-Diagonal
converge quickly with similar speeds, with AdaGrad-Norm being slightly
faster. The remaining algorithms do not converge.

For $\tau=\sigma=1$, the CP algorithm does not converge. Interestingly, AdaGrad-Norm and
AdaGrad-Diagonal both still reach small TV residuals quickly. Adam-Norm
reaches very small TV residuals quickly, but has a very unstable
behavior. Remaining algorithms do not converge.

Overall, the AdaGrad-Norm algorithm demonstrates good robustness and 
better performance than the CP algorithm,
and even achieves a rapid decrease for $\tau=\sigma=1$, where the
CP algorithm is unstable.

\subsection{Solving games with the Mirror-Prox operator}
\label{sec:games-mirror-prox}
Let $m,n\geqslant 1$. We consider a two-player zero-sum game, encoded in a matrix
$B\in \mathbb{R}^{m\times n}$ and its duality gap
\[ \delta(a,b):=\max_{a'\in \Delta_m}\left< a', Bb \right>-\min_{b'\in \Delta_n}\left< a,Bb' \right>. \]
The von Neumann minimax
theorem~\citep{vonneumann1928theorie} ensures
the existence of a solution $(a_*,b_*)\in \Delta_m\times \Delta_n$ such that
$\delta(a_*,b_*)=0$, which implies that
\[ \left< a_*, Bb_* \right> =\max_{a\in \Delta_m}\min_{b\in \Delta_n}\left< a, Bb \right> = \min_{b\in \Delta_n}\max_{a\in \Delta_m}\left< a, Bb \right>.  \]

A classical approach for solving such a game is to cast it as a variational
inequality problem, because the above can be equivalently written as 
\[ \min_{(a,b)\in \Delta_m\times \Delta_n}\left< H(a_*,b_*), (a,b)-(a_*,b_*) \right> \geqslant 0, \]
where operator $H$ is defined as
\[ H(a,b)=(-Bb,B^{\top}\!a),\quad (a,b)\in \Delta_m\times \Delta_n.   \]
Key properties of this operator are Lipschitz continuity, and
monotonicity, in the sense that
\[ \left< H(a',b')-H(a,b),(a',b')-(a,b)  \right>\geqslant 0,\quad (a,b),(a',b')\in \Delta_m\times \Delta_n, \]
which in this case is in fact an equality.

An important family of algorithms for solving variational inequalities
with monotone and Lipschitz continuous operators is Mirror-Prox (MP)
\citep{nemirovski2004prox}, which involves the choice of a
regularizer, aka distance-generating function which defines the
underlying geometry followed by the algorithm. For simplicity, we only
consider the Euclidean instance of MP, which can be rewritten as
KM iterations associated with the following operator
\[ F_\gamma:=2\cdot \Pi_{\Delta_m\times \Delta_n}\circ (I-\gamma (H\circ \Pi_{\Delta_m\times \Delta_n}\circ (I-\gamma H)))-I, \]
where $\gamma>0$ is a step-size.
The MP baseline is $(I+F_\gamma)/2$, namely the projected
extragradient update. The sign convention for $H$ above agrees with
the subtractive steps in the code and with the maximizing role of
the first player.

We perform experiments by sampling a game matrix of size $600\times 400$
and rank $30$, by computing $B=UV^{\top}$ where $U\in\mathbb{R}^{600\times30}$ and
$V\in\mathbb{R}^{400\times30}$ are Gaussian arrays. We then use operator $F_\gamma$ on $\mathcal{X}=\Delta_m\times \Delta_n$.
All methods start from uniform strategies and perform $1000$
iterations. After each iteration, the code evaluates
$\delta(a,b)=\max_i(Bb)_i-\min_j(B^{\top}a)_j$, without averaging
the strategies or normalizing the gap by the scale of $B$.

\begin{figure}[p]
  \centering
  \begin{subfigure}[c]{1.0\linewidth}
  \centering
\resizebox{225pt}{!}{\input{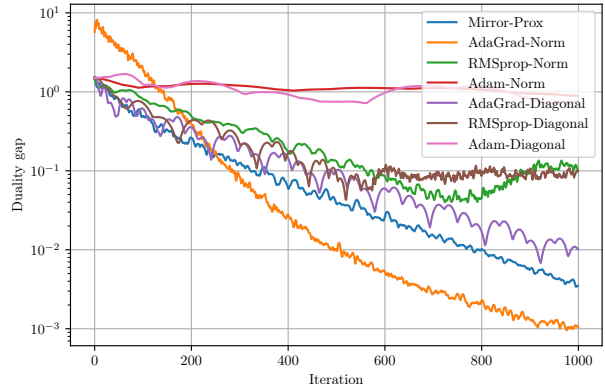}}
\caption{$\gamma=3\cdot 10^{-4}$, the Mirror-Prox algorithm converges slowly.}
  \end{subfigure}
  \par\smallskip
  \begin{subfigure}[c]{1.0\linewidth}
  \centering
\resizebox{225pt}{!}{\input{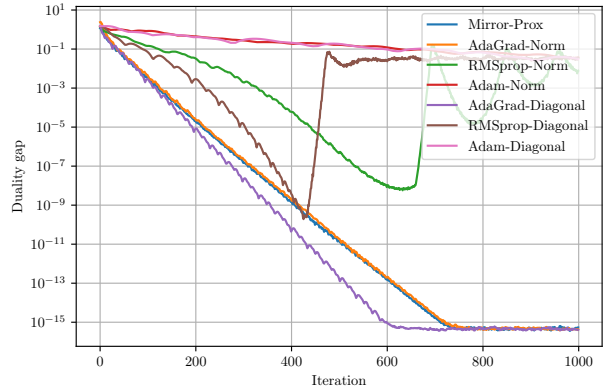}}
\caption{$\gamma=10^{-3}$, the Mirror-Prox algorithm converges quickly.}
  \end{subfigure}
  \par\smallskip
  \begin{subfigure}[c]{1.0\linewidth}
  \centering
\resizebox{225pt}{!}{\input{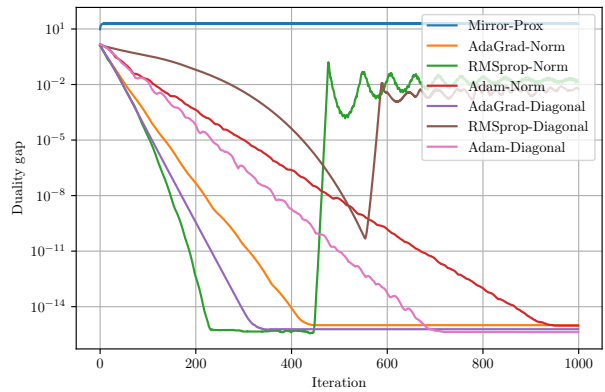}}
\caption{$\gamma=10^{-2}$, the Mirror-Prox algorithm does not converge.}
  \end{subfigure}

\caption{Solving a two-player zero-sum game.}\label{fig:game}
\end{figure}

Simplex projections use bisection with at most $100$ steps and a
$10^{-10}$ tolerance on the total mass. The Diagonal variants use the
square roots of the coordinate-wise accumulators as weights. The
convergence is observed through the duality gap and is plotted in
Figure~\ref{fig:game}. MP is highly sensitive to the choice of the
step-size $\gamma$, for which we try three different values.

For $\gamma=3\cdot 10^{-4}$, the MP algorithm converges slowly. 
Except for initial iterations, AdaGrad-Norm seems to converge slightly
faster than MP.\@ AdaGrad-Diagonal converges slightly more slowly than
MP.\@ The remaining algorithms do not converge.

For $\gamma =10^{-3}$, the MP algorithm converges quickly, and AdaGrad-Norm
achieves a very similar convergence. AdaGrad-Diagonal converges faster
than MP.\@ The remaining algorithms do not converge.

For $\gamma =10^{-2}$, the MP algorithm does not converge. Except for
RMSprop-Diagonal, all other algorithms quickly reach gaps near
floating-point precision. Surprisingly, the fastest algorithm to reach
such small gaps is RMSprop-Norm, which shows instability, as its gap increases
afterwards. AdaGrad-Norm and AdaGrad-Diagonal converge quickly (faster
than MP with $\gamma= 10^{-3}$).

Overall, AdaGrad-Norm and AdaGrad-Diagonal seem to converge similarly
to MP in cases where the latter converges, and still achieve (faster) convergence
in the case where MP does not converge, illustrating robustness
for this instance and these parameter choices.
RMSprop-Norm indicates the possibility of very fast convergence in
the case of a large step-size $\gamma$, but with instability.

\clearpage
\section{Conclusion and perspectives}
\label{sec:conclusion}

We introduced a novel approach for defining and analyzing fixed point
iterations based on regret minimization. This approach recovers the
important Krasnoselskii--Mann iterations as a special case. The great
variety of existing regret minimization algorithms with various
sophisticated regret bounds suggests that this approach could give rise
to many novel fixed-point iterations with interesting convergence
properties, as well as shedding new light on existing ones. We also
hope that this will lead to faster algorithms in practice in various
fields of application.

In this paper, we have obtained simple iterations for non-self
mappings and then 
focused on iterations based on the AdaGrad family of regret minimizers. The conversion of
AdaGrad-Norm into fixed-point iterations
offered convergence guarantees that are \emph{adaptive} to unknown
$L$-nonexpansiveness (meaning to the best-in-hindsight positive scaling).
AdaGrad-Diagonal and AdaGrad-Full achieved adaptivity to
$A$-nonexpansiveness (meaning to the best-in-hindsight positive
diagonal and positive definite
preconditioner, respectively). To the best of our knowledge, these adaptive
guarantees are the first of their kind.

Numerical experiments on three different problems illustrated that
AdaGrad-based iterations can achieve faster and more robust
convergence than baseline methods.

The ideas and results presented in this work suggest a wide array of
research directions---let us mention a few.

\subsection*{Adaptivity to the distance to the solution}
Parameter $\eta>0$ in the AdaGrad-type iterations has to be chosen either
as a function of a known upper bound on the distance to the solution, or
as a function of the diameter of the domain (if finite), in order for
the convergence bound to scale only linearly in the upper bound on the
distance. Without such a known upper bound, the dependence becomes
quadratic. In the context of convex optimization, parameter-free
variants of AdaGrad that adaptively achieve the best dependence on the
distance to the solution have been
proposed~\citep{orabona2018scale,defazio2023learning,ivgi2023dog,khaled2023dowg}.
More generally, such near-linear dependence on the distance to the
solution has been obtained in the parameter-free regret minimization
literature \citep{mcmahan2012no,mcmahan2014unconstrained,orabona2016coin,cutkosky2018black,jacobsen22parameter,cutkosky2024fully}. Adapting such guarantees for fixed points seems to be a
promising research direction.

\subsection*{More AdaGrad variants}
Our numerical experiments include fixed-point iterations converted
from AdaGrad variants RMSprop and Adam: these algorithms are very
successful in deep learning optimization, but with lesser
theoretical understanding. From a practical perspective, we could
further explore the performance of the many other AdaGrad variants \citep{reddi2018convergence,loshchilov2018decoupled,Liu2020On,zhuang2020adabelief} even if they often lack solid theoretical guarantees.

Some other AdaGrad variants provide both solid theoretical guarantees as
well as excellent practical
performance~\citep{kavis2019unixgrad,antonakopoulos2022undergrad}, and
adapting them for fixed points looks very tempting and promising.

\subsection*{Beyond finite dimension}
The present work is restricted to finite dimension, but the theory of
fixed-point iterations is well developed in Hilbert
spaces~\citep{bauschke2011convex} and Banach spaces~\citep{berinde2007iterative}. An important question is whether a
regret-based approach is possible and relevant in infinite dimension.
For AdaGrad-Norm, the residual bounds extend to Hilbert spaces by the
same proof. Under global $L$-nonexpansiveness, the quasi-Fejér argument
gives weak convergence: boundedness provides weak cluster points,
demiclosedness of $I-F_L$ places them in the fixed-point set, and
Opial's lemma yields weak convergence of the whole sequence.
Strong convergence would require additional assumptions.

\subsection*{Adaptivity to contractive properties obtained through scaling}
An important aspect of fixed-point iterations is that contractive
properties lead to geometric convergence as in the Banach--Picard
theorem.
\citet{diakonikolas2025pushing} proposes iterations that guarantee the
optimal complexities both in the contraction and the nonexpansive cases.
An interesting question is whether geometric convergence can
be analyzed through regret minimization, and whether it can be
achieved in settings when an unknown positive (definite) scaling is
necessary to obtain a contractive operator. For instance, operator 
\[ F(x)=Mx, \quad x\in \mathbb{R}^2,\quad \text{where }M=\begin{pmatrix}-\alpha&0\\0&1-\varepsilon
\end{pmatrix},\quad \alpha>2,\quad  0<\varepsilon<2, \]
from Section~\ref{sec:introduction} is not a contraction (not even nonexpansive),
but scaled operator $F_L:=I+\frac{1}{L}(F-I)$ is a $(1+\alpha-\varepsilon)/(1+\alpha+\varepsilon)$-contraction for $L=(1+\alpha+\varepsilon)/2$.
The question is whether there exist iterations that guarantee geometric
convergence in such cases without prior knowledge of $L$.

\subsection*{Adaptive guarantees with $1/T$ convergence}

An important recent breakthrough regarding fixed-point iterations with
nonexpansive operators is the $1/T$ convergence obtained for 
carefully \emph{tuned Halpern} iterations~\citep{halpern1967fixed,lieder2021convergence}. A similar
rate was also achieved using a Nesterov-type momentum applied to KM~\citep{boct2023fast}.
A natural question is whether such fast rates can be recovered using
the regret minimization approach and whether adaptivity to 
unknown positive (definite) scaling can be achieved in addition.

\subsection*{Bundle methods: fixed-point iterations with memory}

In first-order convex optimization, a class of
methods~\citep{lemarechal1975,kiwiel1983aggregate,kiwiel1990proximity}
use (a subset of) previously computed gradients to create a model of
the objective function. For instance, as soon as a gradient of the
objective function $\nabla f(x_t)$ is computed at some point $x_t$, then a
minimizer $x_*$ belongs to the half-space
$\mathcal{X}_t:=\left\{ x\in \mathbb{R}^d,\ \left< x-x_t, \nabla f(x_t) \right> \leqslant 0 \right\}$
because it holds by convexity that
\[ 0\geqslant f(x_*)-f(x_t)\geqslant \left< \nabla f(x_t), x_*-x_t \right>. \] Then by
keeping each previously computed gradient in memory, an algorithm can
speed up convergence by performing at each step a projection onto the
intersection of previous sets $(\mathcal{X}_s)_{1\leqslant s\leqslant t}$.

A similar remark holds with nonexpansive operators $F$ with fixed
point $x_*$. As soon as $F(x_t)$ is computed at some point $x_t$,
solution $x_*$ necessarily belongs to half-space
\[ \mathcal{X}_t:=\left\{ x\in \mathbb{R}^d,\ \left< x, F(x_t)-x_t \right> \geqslant \left< x_t, F(x_t)-x_t \right> + \frac{\left\| x_t-F(x_t) \right\|_2^2}{2} \right\}, \]
as a consequence of Lemma~\ref{lm:online-to-fp}.
Such half-spaces already underpin fixed-point projection methods:
\citet{bauschke2012projection} project the initial point onto their
successive intersections for quasi-nonexpansive mappings.
Related separating half-spaces are used for fallback projections in
SuperMann~\citep{themelis2019supermann}.

Just like projected fixed-point iterates on a constant domain
are naturally analyzed with Projected Online Gradient
Descent in Theorem~\ref{thm:proj-KM}, a natural direction is to analyze
fixed-point iterations on nested shrinking domains through regret minimization.




\section*{Acknowledgments}

The author is indebted to Roberto Cominetti for introducing him to the
topic of fixed point iterations, and Krasnoselskii--Mann iterations in
particular, during the author's stay at Universidad de Chile,
Santiago, in November 2016, funded by ECOS-Sud project No.~\texttt{C15E03}. This work was supported by the French Agence
Nationale de la Recherche (ANR) under reference
\texttt{ANR-21-CE40-0020} (CONVERGENCE project), and by a public grant
as part of the \emph{Investissement d'avenir project}, reference
\texttt{ANR-11-LABX-0056-LMH}, LabEx LMH. This work also benefited
from insightful discussions with Patrick L.\@ Combettes, Sylvain
Sorin, Kfir Levy, Makiel Riveline, Yura Malitsky and from the
in-depth reports by the anonymous referees.


\bibliographystyle{abbrvnat}
\bibliography{bib}

\appendix

\section{Mirror Descent with time-dependent metrics}
\label{sec:mirror-descent-type}
We first state a generic regret bound for Online Mirror Descent
associated with nondecreasing squared Mahalanobis norms that contain as special cases
Projected Online Gradient Descent, AdaGrad-Norm, AdaGrad-Diagonal and AdaGrad-Full.
Specific bounds for each of these special cases are then deduced.

\begin{proposition}[Regret bound for Online Mirror Descent with
  nondecreasing squared Mahalanobis norms]
  \label{prop:omd}
  Let $(u_t)_{t\geqslant 1}$ be a sequence in $\mathbb{R}^d$, $(A_t)_{t\geqslant 0}$ a sequence in $\mathcal{S}_{++}^d(\mathbb{R})$ such that
  $A_{t+1}-A_t$ is positive semidefinite for all $t\geqslant 0$.
  Let $x_1\in \mathcal{X}$ and
  \[ x_{t+1}=\Pi_{\mathcal{X},A_t}(x_t+A_t^{-1}u_t),\quad t\geqslant 1.  \]
  Let $x\in \mathcal{X}$ and $T\geqslant 1$.
\begin{enumerate}[(i)]
\item\label{item:omd} Then, 
\[ \sum_{t=1}^T\left< u_t, x-x_t \right> \leqslant \frac{\left\| x-x_1 \right\|_{A_0}^2}{2}+\frac{1}{2}\sum_{t=1}^T\left( \left\| x-x_{t} \right\|_{A_{t}-A_{t-1}}^2+\left\| u_t \right\|_{A_t^{-1}}^2 \right). \]
\item\label{item:ogd-step-sizes} Let $(\eta_t)_{t\geqslant 1}$ be a positive and
  nonincreasing sequence. If $A_t=\eta_t^{-1}I_d$ for all $t\geqslant 1$, then
  \[ x_{t+1}=\Pi_{\mathcal{X}}\left( x_t+\eta_tu_t \right),\quad t\geqslant 1,  \]
  and
  \[ \sum_{t=1}^T\left< u_t, x-x_t \right> \leqslant \frac{\displaystyle \max_{1\leqslant t\leqslant T}\left\| x-x_t \right\|_2^2}{2\eta_T}+\sum_{t=1}^T\frac{\eta_t\left\| u_t \right\|_2^2}{2}. \]
\end{enumerate}
\end{proposition}
\begin{proof}
Let $t\geqslant 1$. Because $x\in \mathcal{X}$ and the projection operator $\Pi_{\mathcal{X},A_t}$ is
nonexpansive for $\left\|\,\cdot\,\right\|_{A_t}$, using the
definition of $x_{t+1}$ yields
\begin{align*}
\left\| x_{t+1}-x \right\|_{A_t}^2&\leqslant \left\| x_t+A_t^{-1}u_t-x \right\|_{A_t}^2=\left\| x_t-x \right\|_{A_t}^2+2\left< u_t, x_t-x \right> + \left\| u_t \right\|_{A_t^{-1}}^2.
\end{align*}
It follows that
\[ 
\left< u_t, x-x_t \right> \leqslant \frac{\left\| x-x_t \right\|_{A_t}^2}{2}-\frac{\left\| x-x_{t+1} \right\|_{A_{t+1}}^2 }{2}+\frac{\left\| x-x_{t+1} \right\|_{A_{t+1}-A_t}^2}{2}+\frac{\left\| u_t \right\|_{A_t^{-1}}^2}{2}. \]

Then summing yields inequality \eqref{item:omd}. We now turn to \eqref{item:ogd-step-sizes}. Injecting
$A_t=\eta_t^{-1}I_d$ into the general case, we obtain
\[ \sum_{t=1}^T\left< u_t, x-x_t \right>\leqslant\frac{\left\| x-x_1 \right\|_2^2}{2\eta_1}+\sum_{t=2}^T\left( \frac{1}{\eta_t}-\frac{1}{\eta_{t-1}} \right)\frac{\left\| x-x_t \right\|_2^2}{2}+\sum_{t=1}^T\frac{\eta_t\left\| u_t \right\|_{2}^2}{2}. \]
Because the sequence $(\eta_t)_{t\geqslant 1}$ is positive and nonincreasing,
$1/\eta_t-1/\eta_{t-1}\geqslant 0$ and we can bound each squared distance
$\left\| x-x_t \right\|_2^2$ by their maximum and write
\[ \sum_{t=1}^T\left< u_t, x-x_t \right>\leqslant\frac{\displaystyle \max_{1\leqslant t\leqslant T}\left\| x-x_t \right\|_2^2}{2\eta_T}+\sum_{t=1}^T\frac{\eta_t\left\| u_t \right\|_{2}^2}{2}. \]
\end{proof}

\subsection{Projected Online Gradient Descent}
\begin{proposition}[Regret bound for Projected Online Gradient Descent]
  \label{prop:proj-ogd}
 Let \((u_t)_{t\geqslant 1}\) be a sequence in \(\mathbb{R}^d\), \(x_1\in \mathcal{X}\) and
 \begin{equation}
\label{eq:proj-ogd}
x_{t+1}=\Pi_{\mathcal{X}}(x_t+u_t),\quad t\geqslant 1. 
\end{equation}
Then for all integer \(T \geqslant 1\) and $x\in \mathcal{X}$,
 \[ \sum_{t=1}^{T}\left< u_t, x-x_t \right> \leqslant \frac{\left\| x-x_1 \right\|_2^2}{2}+\frac{1}{2}\sum_{t=1}^{T}\left\| u_t \right\|_2^2. \]
\end{proposition}
\begin{proof}
Projected Online Gradient Descent corresponds to Online Mirror
Descent on $\mathcal{X}$ with constant matrix $A_t=I_d$ (for all $t\geqslant 1$).
Then, the iteration reduces to
\[ x_{t+1}=\argmin_{x\in \mathcal{X}}\frac{1}{2}\left\| x-(x_t+u_t) \right\|_2^2=\Pi_{\mathcal{X}}(x_t+u_t),\quad t\geqslant 1, \]
which is exactly Projected Online Gradient Descent.
The regret bound from Proposition~\ref{prop:omd} boils down to
\[ \sum_{t=1}^T\left< u_t, x-x_t \right> \leqslant\frac{\left\| x-x_1 \right\|_2^2}{2}+\sum_{t=1}^T\frac{\left\| u_t \right\|_2^2}{2},  \]
hence the result.
\end{proof}

\subsection{AdaGrad-Norm}
\begin{proposition}[Regret bound for AdaGrad-Norm]
  \label{prop:adagrad-norm-regret-bound}
Let \((u_t)_{t\geqslant 1}\) be a sequence in \(\mathbb{R}^d\), \(\eta>0\), \(x_1\in \mathcal{X}\) and
for \(t\geqslant 1\),
\[ x_{t+1}=\Pi_{\mathcal{X}}\left( x_t+\frac{\eta}{\sqrt{\sum_{s=1}^t\left\| u_s \right\|_2^2}}u_t \right),  \]
with convention $0/0=0$.
Then for all \(x\in \mathcal{X}\) and \(T\geqslant 1\),
\[ \sum_{t=1}^T\left< u_t, x-x_t \right> \leqslant \left( \frac{D_{2,T}^2}{2\eta}+\eta \right)\sqrt{\sum_{t=1}^T\left\| u_t \right\|_2^2}.  \]
where $D_{2,T}=\max_{1\leqslant t\leqslant T}\left\| x_t-x \right\|_2$. In particular,
for \(D\geqslant (\operatorname{diam}_2\mathcal{X})\), choosing $\eta=D/\sqrt{2}$ yields
\[ \sum_{t=1}^T\left< u_t, x-x_t \right> \leqslant D\sqrt{2\sum_{t=1}^T\left\| u_t \right\|_2^2}.  \]
\end{proposition}
\begin{proof}
The proof is adapted from \citet{levy2017online}. If $u_t=0$ for all $t\geqslant 1$, the result holds. Otherwise, let
  \[ \tau=\min_{}\left\{ t\geqslant 1,\ u_t\neq 0 \right\}.  \]
Let $(\eta_t)_{t\geqslant 1}$ be a positive and nonincreasing sequence defined as
\[ \eta_t=
\begin{cases}
\displaystyle \frac{\eta}{\left\| u_{\tau} \right\|_2}&\text{if $t\leqslant \tau$}\\
\displaystyle \frac{\eta}{\sqrt{\sum_{s=1}^t\left\| u_s \right\|_2^2}}&\text{if $t\geqslant \tau$}.  
\end{cases}
 \]
Then, $(x_t)_{t\geqslant 1}$ is a sequence of Projected OGD iterates
on $\mathcal{X}$ with positive and nonincreasing step-sizes
$(\eta_t)_{t\geqslant 1}$ as in Proposition~\ref{prop:omd} \eqref{item:ogd-step-sizes}, because for $1\leqslant t\leqslant \tau-1$, $u_t=0$ and thus, 
\[ x_{t+1}=\Pi_{\mathcal{X}}\left(x_t \right) =\Pi_{\mathcal{X}}\left( x_t+\eta_tu_t \right),  \]
and for $t\geqslant \tau$,
\[ x_{t+1}=\Pi_{\mathcal{X}}\left( x_t+\frac{\eta}{\sqrt{\sum_{s=1}^t\left\| u_s \right\|_2^2}}u_t \right) =\Pi_{\mathcal{X}}\left( x_t+\eta_tu_t \right).  \]
If \(T<\tau\), the result holds because both quantities are zero. If \(T\geqslant \tau\),
the regret bound for Projected OGD with nonincreasing step-sizes from Proposition~\ref{prop:omd} \eqref{item:ogd-step-sizes} gives
\begin{align*}
\sum_{t=1}^T\left< u_t, x-x_t \right> &\leqslant \frac{\displaystyle \max_{1\leqslant t\leqslant T}\left\| x-x_t \right\|_2^2}{2\eta_T}+\sum_{t=1}^T\frac{\eta_t\left\| u_t \right\|_2^2}{2}\\
&=\frac{\displaystyle \max_{1\leqslant t\leqslant T}\left\| x_t-x \right\|_2^2}{2\eta}\sqrt{\sum_{t=1}^T\left\| u_t \right\|_2^2}+\frac{\eta}{2}\sum_{t=1}^T\frac{\left\| u_t \right\|_2^2}{\sqrt{\sum_{s=1}^t\left\|u_s \right\|_2^2}}\\
  &\leqslant \left( \frac{\max_{1\leqslant t\leqslant T}\left\| x_t-x \right\|_2^2}{2\eta}+\eta\right)\sqrt{\sum_{t=1}^T\left\| u_t \right\|_2^2}. 
\end{align*}
using Lemma~\ref{lm:discrete-int-inv-sqrt} below. Hence the result.
\end{proof}
\begin{lemma}[Lemma 3.5 in \citet{auer2002adaptive}]
  \label{lm:discrete-int-inv-sqrt}
Let $(a_t)_{t\geqslant 1}$ be a nonnegative sequence. Then for all $T\geqslant 1$,
\[ \sum_{t=1}^T\frac{a_t}{\sqrt{\sum_{s=1}^ta_s}}\leqslant 2\sqrt{\sum_{t=1}^Ta_t}. \]
with convention $0/0=0$.
\end{lemma}

\subsection{AdaGrad-Diagonal}
\begin{proposition}[Regret bound for AdaGrad-Diagonal]
  \label{prop:adagrad-diagonal}
Let \((u_t)_{t\geqslant 1}\) be a sequence in \(\mathbb{R}^d\), \(\eta,\varepsilon>0\), \(x_1\in \mathcal{X}\)
and for \(t\geqslant 1\),
\[ x_{t+1}=\Pi_{\mathcal{X},A_t}(x_t+A^{-1}_tu_t),  \]
where
\[ A_t=\eta^{-1}\operatorname{diag}\left( \sqrt{\varepsilon^2+\sum_{s=1}^tu_{s,i}^2} \right)_{1\leqslant i\leqslant d}.   \]
Let $x\in \mathcal{X}$ and $T\geqslant 1$.
\begin{enumerate}[(i)]
\item\label{item:adagrad-diagonal-regret} Then,
\begin{multline*}
\sum_{t=1}^T\left< u_t, x-x_t \right>
\leqslant \frac{\varepsilon}{2\eta}\left\| x_1-x \right\|_{2}^2\\
{}+\left( \frac{D_{\infty,T}^2}{2\eta} + \eta\right)
\inf_{\substack{A\in \mathcal{S}_{++}^d(\mathbb{R})\\A\textup{ diagonal}}}
\sqrt{(\operatorname{Tr}A)\sum_{t=1}^T\left\| u_t \right\|_{A^{-1}}^2},
\end{multline*}
where $D_{\infty,T}=\max_{1\leqslant t\leqslant T}\left\| x_t-x \right\|_{\infty}$.
\item\label{item:adagrad-diagonal-radius} Let 
  $D\geqslant \operatorname{diam}_{\infty}(\mathcal{X})$. Then choosing $\eta=D/\sqrt{2}$ yields
\[ \sum_{t=1}^T\left< u_t, x-x_t \right>\leqslant D\left( \frac{\varepsilon d}{\sqrt{2}}+\inf_{\substack{\substack{A\in \mathcal{S}_{++}^d(\mathbb{R})\\A\textup{ diagonal}}}}\sqrt{2(\operatorname{Tr}A)\sum_{t=1}^T\left\| u_t \right\|_{A^{-1}}^2} \right).  \]
\end{enumerate}
\end{proposition}
\begin{proof}
The proof is adapted from \citet{duchi2011adaptive}. Let $A_0=\eta^{-1}\varepsilon I_d$.
  Then, $(x_t)_{t\geqslant 1}$ is the sequence of Online Mirror Descent 
  iterates on $\mathcal{X}$ associated with positive diagonal matrices $(A_t)_{t\geqslant 1}$.
  Then, Proposition~\ref{prop:omd} gives
  \begin{equation}
\label{eq:14}
\sum_{t=1}^T\left< u_t, x-x_t \right> \leqslant \frac{\varepsilon\left\| x-x_1 \right\|_{2}^2}{2\eta}+\frac{1}{2}\sum_{t=1}^T\left\| u_t \right\|_{A_t^{-1}}^2+\frac{1}{2}\sum_{t=1}^T\left\| x-x_{t} \right\|_{A_{t}-A_{t-1}}^2. 
\end{equation}
Because matrices $(A_t)_{t\geqslant 0}$ are diagonal and nondecreasing, the above last sum is bounded as
\begin{align}
  \label{eq:12}
  \begin{split}
\frac{1}{2}\sum_{t=1}^T\left\| x-x_t \right\|_{A_t-A_{t-1}}^2&=\frac{1}{2}\sum_{t=1}^T\sum_{i=1}^d(x_i-x_{t,i})^2(A_{t,ii}-A_{t-1,ii})\\
  &\leqslant \frac{\displaystyle \max_{1\leqslant t\leqslant T}\left\| x-x_t \right\|_{\infty}^2}{2}\sum_{i=1}^d \sum_{t=1}^T(A_{t,ii}-A_{t-1,ii})\\
  &=\frac{\displaystyle \max_{1\leqslant t\leqslant T}\left\| x-x_t \right\|_{\infty}^2}{2}\sum_{i=1}^d(A_{T,ii}-A_{0,ii})\\
  &=\frac{\displaystyle \max_{1\leqslant t\leqslant T}\left\| x-x_t \right\|_{\infty}^2}{2\eta}\sum_{i=1}^d\left( \sqrt{\varepsilon^2+\sum_{t=1}^Tu_{t,i}^2}-\varepsilon \right) \\
  &\leqslant \frac{\displaystyle \max_{1\leqslant t\leqslant T}\left\| x-x_t \right\|_{\infty}^2}{2\eta}\sum_{i=1}^d\sqrt{\sum_{t=1}^Tu_{t,i}^2}  .
  \end{split}
\end{align}
Regarding the second term from the right-hand side of \eqref{eq:14},
\begin{align}
\label{eq:10}
  \begin{split}
\frac{1}{2}\sum_{t=1}^T\left\| u_t \right\|_{A_t^{-1}}^2&=\frac{\eta}{2}\sum_{t=1}^T\sum_{i=1}^d\frac{u_{t,i}^2}{\sqrt{\varepsilon^2+\sum_{s=1}^tu_{s,i}^2}}\\
  &\leqslant \frac{\eta}{2}\sum_{t=1}^T\sum_{i=1}^d\frac{u_{t,i}^2}{\sqrt{\sum_{s=1}^tu_{s,i}^2}}\\
    &\leqslant \eta\sum_{i=1}^d\sqrt{\sum_{t=1}^Tu_{t,i}^2}\\
    &\leqslant \eta \inf_{\substack{A\in \mathcal{S}_{++}^d(\mathbb{R})\\\text{$A$ diagonal}}}\sqrt{(\operatorname{Tr}A)\sum_{t=1}^T\left\| u_t \right\|_{A^{-1}}^2}, 
  \end{split}
\end{align}
where we used Lemma~\ref{lm:discrete-int-inv-sqrt}  (resp.\
Lemma~\ref{lm:sum-sqrt} below) for the penultimate inequality (resp.\ the last inequality). Then combining \eqref{eq:14},
\eqref{eq:12} and \eqref{eq:10} gives the result.
\end{proof}

\begin{lemma}
  \label{lm:sum-sqrt}
  Let $b=(b_1,\dots,b_d)\in \mathbb{R}_+^d$. Then,
  \[ \sum_{i=1}^db_i=\inf_{\substack{A\in \mathcal{S}_{++}^d(\mathbb{R})\\\text{$A$ diagonal}}}\sqrt{(\operatorname{Tr}A)\left\| b \right\|_{A^{-1}}^2}. \]
\end{lemma}
\begin{proof}
Let $a_1,\dots,a_d>0$. Then using Jensen's
inequality for $z\mapsto z^2$,
\[ \left( \sum_{i=1}^d\left( \frac{a_i}{\sum_{j=1}^da_j} \right)\left( \frac{b_i}{a_i} \right)   \right)^2 \leqslant \sum_{i=1}^d\frac{a_i}{\sum_{j=1}^da_j}\frac{b_i^2}{a_i^2}=\sum_{i=1}^d\frac{b_i^2/a_i}{\sum_{j=1}^da_j}. \]
Taking the square root gives
\[ \sum_{i=1}^db_i\leqslant \sqrt{\left( \sum_{j=1}^da_j \right) \sum_{i=1}^d\frac{b_i^2}{a_i}}=\sqrt{(\operatorname{Tr}A)\left\| b \right\|_{A^{-1}}^2}, \]
where $A=\operatorname{diag}(a_1,\dots,a_d)$ is a positive diagonal
matrix.  This proves
  \[ \sum_{i=1}^db_i\leqslant \inf_{\substack{A\in \mathcal{S}_{++}^d(\mathbb{R})\\\text{$A$ diagonal}}}\sqrt{(\operatorname{Tr}A)\left\| b \right\|_{A^{-1}}^2}. \]
Let $\varepsilon>0$ and consider
\[ a_i=b_i+\varepsilon,\quad 1\leqslant i\leqslant d, \]
which are always positive. Then,
\[ \sqrt{(\operatorname{Tr}A)\sum_{i=1}^d\frac{b_i^2}{a_i}}=\sqrt{\left( \sum_{i=1}^db_i+d\varepsilon \right)\sum_{i=1}^d\frac{b_i^2}{b_i+\varepsilon} } \]
which converges to $\sum_{i=1}^db_i$ as $\varepsilon\to 0^{+}$. This proves
  \[ \sum_{i=1}^db_i\geqslant \inf_{\substack{A\in \mathcal{S}_{++}^d(\mathbb{R})\\\text{$A$ diagonal}}}\sqrt{(\operatorname{Tr}A)\left\| b \right\|_{A^{-1}}^2}. \]
Hence the result.
\end{proof}

\subsection{AdaGrad-Full}
\begin{proposition}[Regret bound for AdaGrad-Full]
\label{prop:adagrad-full}
Let \((u_t)_{t\geqslant 1}\) be a sequence in \(\mathbb{R}^d\), \(\eta,\varepsilon>0\), \(x_1\in \mathcal{X}\)
and for \(t\geqslant 1\),
\[ x_{t+1}=\Pi_{\mathcal{X},A_t}(x_t+A^{-1}_tu_t),  \]
where
\[ A_t=\eta^{-1}\left( \varepsilon^2 I_d+\sum_{s=1}^tu_su_s^{\top} \right)^{1/2}.   \]
Let $x\in \mathcal{X}$ and $T\geqslant 1$.
\begin{enumerate}[(i)]
\item\label{item:adagrad-full-regret} Then,
  \begin{multline*}
    \sum_{t=1}^T\left< u_t, x-x_t \right> \leqslant \frac{\varepsilon}{2\eta}\left\| x_1-x \right\|_{2}^2+\left( \frac{\max_{1\leqslant t\leqslant T}\left\| x_t-x \right\|_{2}^2}{2\eta} + \eta\right)\\
    \times \inf_{\substack{\substack{A\in \mathcal{S}_{++}^d(\mathbb{R})}}}\sqrt{(\operatorname{Tr}A)\sum_{t=1}^T\left\| u_t \right\|_{A^{-1}}^2}.
  \end{multline*}
\item\label{item:adagrad-full-radius} Moreover, let
  $D\geqslant \operatorname{diam}_2(\mathcal{X})$. Then choosing
$\eta=D/\sqrt{2}$ yields,
\[ \sum_{t=1}^T\left< u_t, x-x_t \right>\leqslant D\left( \frac{\varepsilon }{\sqrt{2}}+\inf_{\substack{\substack{A\in \mathcal{S}_{++}^d(\mathbb{R})}}}\sqrt{2(\operatorname{Tr}A)\sum_{t=1}^T\left\| u_t \right\|_{A^{-1}}^2} \right).  \]
\end{enumerate}
\end{proposition}
\begin{proof}

The proof is adapted from \citet{duchi2011adaptive}. Similarly to the proof of the regret bound for AdaGrad-Diagonal,
applying Proposition~\ref{prop:omd} with 
$A_0=\eta^{-1}\varepsilon I_d$ yields
\begin{equation}
\label{eq:5}
\sum_{t=1}^T\left< u_t, x-x_t \right>  \leqslant \frac{\varepsilon\left\| x-x_1 \right\|_{2}^2}{2\eta}+\frac{1}{2}\sum_{t=1}^T\left\| x-x_t \right\|_{A_t-A_{t-1}}^2+\frac{1}{2}\sum_{t=1}^T\left\| u_t \right\|_{A_t^{-1}}^2. 
\end{equation}
Because the matrix square root is monotone in the Loewner order, it holds
that $A_t-A_{t-1}\succeq 0$,
and above second term can thus be bounded from above as
\begin{align*}
\frac{1}{2}\sum_{t=1}^T\left\| x-x_t \right\|_{A_t-A_{t-1}}^2&=\frac{1}{2}\sum_{t=1}^T\left< x-x_t, (A_t-A_{t-1})(x-x_t) \right> \\
  &\leqslant \frac{1}{2}\sum_{t=1}^T\left\| x-x_t \right\|_2^2\cdot \lambda_{\text{max}}(A_t-A_{t-1})\\
  &\leqslant \frac{1}{2}\sum_{t=1}^T\left\| x-x_t \right\|_2^2\cdot \operatorname{Tr}(A_t-A_{t-1})\\
  &\leqslant \frac{1}{2}\max_{1\leqslant t\leqslant T}\left\| x-x_t \right\|_2^2\cdot \operatorname{Tr}(A_T-A_0)\\
  &\leqslant \frac{\displaystyle \max_{1\leqslant t\leqslant T}\left\| x-x_t \right\|_2^2}{2\eta}\cdot \operatorname{Tr}\left(\left(\sum_{t=1}^Tu_tu_t^{\top}  \right)^{1/2} \right).
\end{align*} 

Using Lemma~\ref{lm:discrete-int-full-matrices} below, the last term
from \eqref{eq:5} is bounded as
\begin{align*}
\frac{1}{2}\sum_{t=1}^T\left\| u_t \right\|_{A_t^{-1}}^2&=\frac{1}{2}\sum_{t=1}^T\operatorname{Tr}(A_t^{-1}u_tu_t^{\top})\\
  &=\frac{\eta^2}{2}\sum_{t=1}^T\operatorname{Tr}(A_t^{-1}(A_t^2-A_{t-1}^2))\\
  &\leqslant \eta^2\sum_{t=1}^T\operatorname{Tr}(A_t-A_{t-1})\\
  &=\eta^2\left( \operatorname{Tr}(A_T)-\operatorname{Tr}(A_0) \right)\\
  &=\eta\operatorname{Tr}\left( \left( \varepsilon^2I_d+\sum_{t=1}^Tu_tu_t^{\top} \right)^{1/2}-\varepsilon I_d  \right).
\end{align*}
Using the spectral decomposition of positive semidefinite matrix
$\sum_{t=1}^Tu_tu_t^{\top}=P\operatorname{diag}(\lambda_1,\dots,\lambda_d)P^{\top}$, where
$P$ is an orthogonal matrix, the above last quantity is equal to
\begin{align*}
&\eta \operatorname{Tr}\left( \left( \varepsilon^2I_d+\operatorname{diag}(\lambda_1,\dots,\lambda_d) \right)^{1/2}-\varepsilon I_d  \right)\\
&=\eta \operatorname{Tr}\left( \operatorname{diag}\left(\sqrt{\varepsilon^2+\lambda_i}-\varepsilon\right)_{1\leqslant i\leqslant d} \right)\\
&\leqslant \eta \operatorname{Tr}\left( \operatorname{diag}\left(\sqrt{\lambda_1},\dots,\sqrt{\lambda_d}\right) \right)\\
  &=\eta \operatorname{Tr}\left( \left( \sum_{t=1}^Tu_tu_t^{\top} \right)^{1/2}  \right). 
\end{align*}

Besides, applying Lemma~\ref{lm:inf-full-matrices} below with
$M=\sum_{t=1}^Tu_tu_t^{\top}$ ensures that
\[ \operatorname{Tr}\left( \left(\sum_{t=1}^Tu_tu_t^{\top}  \right)^{1/2}  \right)=\inf_{A\in \mathcal{S}_{++}^d(\mathbb{R})}\sqrt{(\operatorname{Tr}A)\sum_{t=1}^T\left\| u_t \right\|_{A^{-1}}^2} .  \]
because then for all $A\in \mathcal{S}_{++}^d(\mathbb{R})$, quantity $\operatorname{Tr}(A^{-1}M)$ that appears in the
lemma rewrites as
\begin{align*}
\operatorname{Tr}(A^{-1}M)&=\operatorname{Tr}\left( A^{-1}\sum_{t=1}^Tu_tu_t^{\top} \right)=\sum_{t=1}^T\operatorname{Tr}\left( A^{-1}u_tu_t^{\top} \right)\\
  &=\sum_{t=1}^T\operatorname{Tr}(u_t^{\top}\!A^{-1}u_t)=\sum_{t=1}^Tu_t^{\top}\!A^{-1}u_t\\
  &=\sum_{t=1}^T\left\| u_t \right\|_{A^{-1}}^2.    
\end{align*}
Injecting the above inequalities into \eqref{eq:5} gives the result.
\end{proof}

\begin{lemma}
  \label{lm:discrete-int-full-matrices}
  Let $A,B\in \mathcal{S}_{++}^d(\mathbb{R})$. Then,
  \[ \operatorname{Tr}(A^{1/2})\leqslant \operatorname{Tr}(B^{1/2})+\frac{1}{2}\operatorname{Tr}(B^{-1/2}(A-B)). \]
\end{lemma}
\begin{proof}
The inequality is a special case of Klein's trace inequality
\citep[Theorem 2.11]{carlen2010trace} applied with concave real-valued
function $x\mapsto x^{1/2}$.
\end{proof}

\begin{lemma}
\label{lm:inf-full-matrices}  
Let $M$ be a real symmetric positive semidefinite matrix of size $d$.
Then,
\[ \operatorname{Tr}\left( M^{1/2} \right)=\inf_{A\in \mathcal{S}_{++}^d(\mathbb{R})}\sqrt{(\operatorname{Tr}A)\operatorname{Tr}(A^{-1}M)}.  \]
\end{lemma}
\begin{proof}
  Let $A\in \mathcal{S}_{++}^d(\mathbb{R})$. Applying the Cauchy--Schwarz inequality for the
  Frobenius inner product,
\begin{align*}
\operatorname{Tr}(M^{1/2})=\operatorname{Tr}(A^{1/2}(A^{-1/2}M^{1/2}))\leqslant \sqrt{\operatorname{Tr}(A)\operatorname{Tr}(A^{-1}M)},
\end{align*}
which proves
\[ \operatorname{Tr}\left( M^{1/2} \right)\leqslant \inf_{A\in \mathcal{S}_{++}^d(\mathbb{R})}\sqrt{(\operatorname{Tr}A)\operatorname{Tr}(A^{-1}M)}.  \]

Conversely, let $\varepsilon>0$, $r$ be the rank of $M$, $\lambda_1(M),\dots,\lambda_d(M)$
be the eigenvalues of $M$ with multiplicity in nonincreasing order,
and $P$ be an orthogonal matrix such that
\[ P^{\top}\!MP=\begin{pmatrix}\lambda_1(M)&&0\\
&\ddots &\\
0&&\lambda_d(M)\end{pmatrix}. \]
Then consider positive definite matrix
\[ A=P\operatorname{diag}(\sqrt{\lambda_1(M)},\dots,\sqrt{\lambda_r(M)},\varepsilon,\dots,\varepsilon)P^{\top}, \]
which implies that
\[ A^{-1}M=P\operatorname{diag}(\sqrt{\lambda_1(M)},\dots,\sqrt{\lambda_r(M)},0,\dots,0)P^{\top}. \]
Hence,
\[ \sqrt{(\operatorname{Tr}A)(\operatorname{Tr}(A^{-1}M))}=\sqrt{(\operatorname{Tr}M^{1/2}+(d-r)\varepsilon)(\operatorname{Tr}M^{1/2})}, \]
which converges to $\operatorname{Tr}M^{1/2}$ as $\varepsilon\to 0^+$ and this
proves
\[ \operatorname{Tr}\left( M^{1/2} \right)\geqslant \inf_{A\in \mathcal{S}_{++}^d(\mathbb{R})}\sqrt{(\operatorname{Tr}A)\operatorname{Tr}(A^{-1}M)}.  \]
Hence the result.
\end{proof}

\section{Extension to quasi-nonexpansiveness}
\label{sec:extens-quasi-nonexp}

We here state the extensions of the results from
Sections~\ref{sec:basic-prop-oper} and \ref{sec:regr-minim-key} to
quasi-nonexpansiveness. The basic notion of quasi-nonexpansiveness
first appears in \citep{dotson1972fixed}. We extend it as follows.

\begin{definition}[Quasi-nonexpansiveness]
Let  \(F:\mathcal{X}\to \mathbb{R}^d\), $x_*\in \mathcal{X}$ a fixed point of $F$ and \(\left\|\,\cdot\,\right\| \) a norm in
\(\mathbb{R}^d\).
\begin{enumerate}[(i)]
\item Operator \(F\) is \emph{quasi-nonexpansive} with respect to
  $x_*$ and \(\left\|\,\cdot\,\right\|\) if for all
\(x\in \mathcal{X}\),
\[ \left\| F(x)-F(x_*) \right\| \leqslant \left\| x-x_* \right\|. \]
\item Let $L>0$. Operator $F$ is \emph{$L$-quasi-nonexpansive} with
  respect to $x_*$ if
  $I+\frac{1}{L}(F-I)$ is quasi-nonexpansive with respect to $x_*$
  and $\left\|\,\cdot\,\right\|_2$.
\item Let $A\in \mathcal{S}_{++}^d(\mathbb{R})$. Operator $F$ is
  \emph{$A$-quasi-nonexpansive} with respect to $x_*$ if
  $I+A^{-1}(F-I)$ is quasi-nonexpansive with respect to $x_*$ and $\left\|\,\cdot\,\right\|_A$.
\end{enumerate}
\end{definition}

We extend to arbitrary norms the notion of star-co-coercivity
introduced by \citep{gorbunov2022extragradient}.
\begin{definition}[Star-co-coercivity]
Let \(G:\mathcal{X}\to \mathbb{R}^d\), $x_*\in \mathcal{X}$ a zero of $G$, \(\left\|\,\cdot\,\right\|\) a norm in \(\mathbb{R}^d\), and \(L>0\).
  Operator \(G\) is \(L\)-\emph{star-co-coercive} with respect to
  $x_*$ and \(\left\|\,\cdot\,\right\|\) if for all \(x\in \mathcal{X}\),
  \[ \left< G(x), x-x_* \right> \geqslant \frac{1}{L}\left\| G(x) \right\|^2_*, \]
  where $\left\|\,\cdot\,\right\|_*$ denotes the dual norm of $\left\|\,\cdot\,\right\|$.
\end{definition}

The following characterization is an extension of \citet[Lemmas C.5 and C.6]{gorbunov2022extragradient}.
\begin{proposition}
  \label{prop:characterization-quasi-nonexp}
  Let \(F:\mathcal{X}\to \mathbb{R}^d\), $x_*\in \mathcal{X}$ a fixed point of $F$, \(G=\frac{I-F}{2}\), and
  \(A\in \mathcal{S}_{++}^d(\mathbb{R})\).
  Then, $F$ is
  $A$-quasi-nonexpansive with respect to $x_*$ if, and only if, \(G\)
  is $1$-star-co-coercive with respect to $x_*$ and
  $\left\|\,\cdot\,\right\|_A$, in other words,
\begin{equation}
\forall x\in \mathcal{X},\quad \left< G(x), x-x_* \right>\geqslant \left\| G(x) \right\|^2_{A^{-1}}. 
\end{equation}
\end{proposition}
\begin{proof}
  Note that \(I+A^{-1}(F-I)=I-2A^{-1}G\). This operator being
  quasi-nonexpansive with respect to $x_*$ and \(\left\|\,\cdot\,\right\|_A\) writes
  \[ \left\| (I-2A^{-1}G)(x)-(I-2A^{-1}G)(x_*) \right\|_A^2\leqslant \left\| x-x_* \right\|_A^2,\quad x\in \mathcal{X}. \]
  Because $G(x_*)=0$ by definition of $x_*$, the above rewrites as
  \begin{equation}
\label{eq:6}
\left\| (I-2A^{-1}G)(x)-x_* \right\|_A^2\leqslant \left\| x-x_* \right\|_A^2,\quad x\in \mathcal{X}. 
\end{equation}
Let \(x\in \mathcal{X}\) be given. Expanding the above left-hand side, we obtain
\begin{multline*}
\left\| (I-2A^{-1}G)(x)-x_* \right\|_A^2=4\left\| A^{-1}G(x) \right\|_A^2\\
\qquad -4\left< A^{-1}G(x), A(x-x_*) \right> + \left\| x-x_* \right\|_A^2\\
=4\left\| G(x) \right\|^2_{A^{-1}}-4\left< G(x), x-x_* \right>+\left\| x-x_* \right\|_A^2. 
\end{multline*}
Plugging the above identity into inequality \eqref{eq:6} and rearranging gives the result.
\end{proof}

\begin{lemma}
  \label{lm:online-to-fp-quasi-nonexp}
  Let $A\in \mathcal{S}_{++}^d(\mathbb{R})$,
  \(F:\mathcal{X}\to \mathbb{R}^d\), \(x_*\in \mathcal{X}\) a fixed point of \(F\). If $F$ is
  $A$-quasi-nonexpansive with respect to $x_*$, then for all \(x\in \mathcal{X}\),
  \[ \left\| F(x)-x \right\|_{A^{-1}}^2\leqslant 2\left< F(x)-x, x_*-x \right>. \]
\end{lemma}
\begin{proof}
Apply Proposition~\ref{prop:characterization-quasi-nonexp} and argue as in the proof of Lemma~\ref{lm:online-to-fp}.
\end{proof}


\end{document}